\documentclass[11pt]{article}
\usepackage[T1]{fontenc}
\usepackage[utf8]{inputenc}
\usepackage[cyr]{aeguill}
\usepackage[english]{babel}
\usepackage{amsmath,amssymb,amsthm}
\usepackage{bm}
\usepackage{graphicx}

\usepackage{authblk}	
\usepackage{geometry}
\usepackage{mathrsfs}

\usepackage{multirow}
\usepackage[justification=centering]{caption}
\usepackage{subfig}
\usepackage{xcolor,colortbl}
\usepackage{diffcoeff}

\usepackage[colorlinks=true,linkcolor=black, citecolor=blue, urlcolor=blue]{hyperref}

\DeclareMathOperator{\Tr}{Tr}

\newtheorem{theorem}{Theorem}
\newtheorem{remark}{Remark}

\newcommand{\sqboxs}{1.2ex}
\newcommand{\sqbox}[1]{\textcolor{#1}{\rule{\sqboxs}{\sqboxs}}}

\definecolor{Blue}{rgb}{0,1,1}
\definecolor{Green}{rgb}{0.4,1,0.4}
\definecolor{GreenYellow}{rgb}{0.3,0.6,0.5}
\definecolor{Yellow}{rgb}{1,1,0.4}
\definecolor{Orange}{rgb}{1,0.6,0.2}
\definecolor{Purple}{rgb}{1,0.5,0.5}
\definecolor{Red}{rgb}{1,0.2,0.2}

\definecolor{Grey}{rgb}{0.5,0.5,0.5}

\newcolumntype{g}{>{\columncolor{Green}[.95\tabcolsep]}c}
\newcolumntype{m}{>{\columncolor{Yellow}[.95\tabcolsep]}c}
\newcolumntype{b}{>{\columncolor{Blue}[.95\tabcolsep]}c}

\graphicspath{{./Figures/}, {./Figures/Eig_Kirchh/}, {./Figures/Simulations/}}

\def\onedot{$\mathsurround0pt\ldotp$}
\def\cddot{
	\mathbin{\vcenter{\baselineskip.67ex
			\hbox{\onedot}\hbox{\onedot}}%
}}

\makeatletter \renewcommand\d[1]{\ensuremath{%
  \;\mathrm{d}#1\@ifnextchar\d{\!}{}}}
\makeatother

\title{Port-Hamiltonian formulation and \\ symplectic discretization of plate models\\
	Part II : Kirchhoff model for thin plates}	

\author[1]{Andrea Brugnoli\thanks{andrea.brugnoli@isae.fr}}
\author[1]{Daniel Alazard\thanks{daniel.alazard@isae.fr}}
\author[1]{Val\'erie Pommier-Budinger \thanks{valerie.budinger@isae.fr}}
\author[1]{Denis Matignon \thanks{denis.matignon@isae.fr}}
\affil[1]{ISAE-SUPAERO, Universit\'e de Toulouse, France. \\
	10 Avenue Edouard Belin, BP-54032, 31055 Toulouse Cedex 4.}

\begin{document}
\maketitle
\begin{abstract}
	The mechanical model of a thin plate with boundary control and 
	observation is presented as a port-Hamiltonian system (PHs\footnote{PHs stands for port-Hamiltonian systems.}), both in vectorial and tensorial forms: the Kirchhoff-Love model of a plate is described by using a Stokes-Dirac structure and this represents a novelty with respect to the existing literature. This formulation is carried out both in vectorial and tensorial forms. Thanks to tensorial calculus, this model is found to mimic the interconnection structure of its one-dimensional counterpart, i.e. the Euler-Bernoulli beam. \\
	The Partitioned Finite Element Method (PFEM\footnote{PFEM stands for partitioned finite element method.}) is then extended to obtain a suitable, i.e. structure-preserving, weak form. The discretization procedure, performed on the vectorial formulation, leads to a finite-dimensional port-Hamiltonian system.  This part II of the companion paper extends part I, dedicated to the Mindlin model for thick plates. The thin plate model comes along with additional difficulties, because of the higher order of the differential operator under consideration.
\end{abstract}

\section*{Introduction}
As presented in part I of this companion paper, the port-Hamiltonian (PH) formalism \cite{BookZwart, bookPHs, Villegas} allows the structured modeling and discretization of multi-physics applications involving interconnected finite- and infinite-dimensional systems \cite{Cervera2007,ShaftIntInfinite}. Preserving the port-Hamiltonian structure in the discretization process is a keypoint to take benefit of this powerful formalism. This issue was first addressed in \cite{Golo}, with a mixed finite element spatial discretization method, and in \cite{moulla:hal-01625008}, with pseudo-spectral methods relying on higher-order global polynomial approximations. All those methods are difficult to implement, especially for those system the spatial dimension of which is bigger than one. Very recently weak formulations which lead to Galerkin numerical approximations began to be explored: in \cite{WeakForm_Kot}, a structure preserving finite element method was introduced for the wave equation in a two-dimensional domain; this method exhibits good results, both in the spectral analysis and simulation part, though requiring of a primal and a dual mesh on the geometry of the problem. Another approach is the partitioned finite element method (PFEM) proposed in \cite{CardosoRibeiro2018}, already largely explored in part I of this companion paper. The advantages of this latter methodology are its simplicity of implementation and its potential to carry over to a wide set of examples, no matter the spatial dimension of the problem. The possible use of open source software like FEniCS \cite{LoggMardalEtAl2012} {or Firedrake \cite{firedrake}} is also an appealing feature of this latter method  \\

In part II of this companion paper, the modeling and discretization of thin plates described by the Kirchhoff-Love plate model is carried out within the PH framework, allowing for boundary control and observation. {The existing literature dealing with the symplectic Hamiltonian formulation of the Kirchhoff plate \cite{LI2016984,LI2018310} focused mainly on analytical solution for the free vibration problem. This approach is powerful whenever easy solution are sought for but does not extend to systems interconnected in complex manners.} Furthermore, plate models were investigated withing the port-Hamiltonian framework using jet theory  {\cite{jetMin,jetKirchh}}, but the numerical implementation of such models remains cumbersome. The main contribution of this paper concerns the representation of the Kirchhoff plate using the concept of Stokes-Dirac structure, so to take advantage from the modularity of this geometric structure. This formalism is presented both in vectorial and tensorial forms.  Moreover, the tensorial formalism \cite[Chapter~16]{Grinfield} highlights that this model mimics the interconnection structure of its one-dimensional counterpart, i.e. the Euler-Bernoulli beam. Compared to part I dedicated to thick plate Mindlin model in which first-order differential operators are explored in dimension two, and compared to \cite{LeGorrec2005} in which second-  or higher-order differential operators were explored in dimension one only, the contribution of this paper is the PH formalism of systems of dimension two described with second-order  differential operators, such as the Kirchhoff-Love  model. The model, once written in a tensorial form, highlights new insights on second-order differential operators: especially the double divergence and the Hessian are proved to be adjoint operators one of another, which represents another important contribution of this paper.
Finally, the extension of the PFEM method  to the structure-preserving discretization of the Kirchhoff model  is also a novelty of the paper. It allows simple implementation of numerical schemes compared to the jet theory formalism, while preserving the structure of PHS at the discrete level. {The last section is dedicated to numerical studies of this model using Firedrake \cite{firedrake}.}

\section{Second-order distributed PH systems: Euler-Bernoulli beam}	
The Euler-Bernoulli beam is the one-dimensional equivalent of the Kirchhoff-Love plate.  This model consists of one PDE, describing the vertical displacement along the beam length:
\begin{equation}
\rho(x) \diffp[2]{w}{t}(x,t) + \displaystyle \diffp[2]{}{x} \left( EI(x) \diffp[2]{w}{x} \right) = 0, \quad x \in (0,L),\, t \ge 0 
\end{equation}
where ${w}(x,t)$ is the transverse displacement of the beam. The coefficients $\rho(x), E(x)$ and $I(x)$  are the mass per unit length, Young's modulus of elasticity and the moment of inertia of a cross section. The energy variables are then chosen as follows:
\begin{equation}
\begin{aligned}
\alpha_{w} &= \rho(x) \diffp{w}{t}(x,t),  &\quad \text{Linear Momentum},\\
\alpha_{\kappa} &= \diffp[2]{w}{x}(x,t), &\quad \text{Curvature}. \\
\end{aligned}
\end{equation}

Those variables are collected in the vector $\bm{\alpha} = (\alpha_{w}, \, \alpha_{\kappa})^T $, so that the Hamiltonian can be written as a quadratic functional in the energy variables: 
\begin{equation}
H = \frac{1}{2} \int_{0}^{L} \bm{\alpha}^T Q \bm{\alpha} \; \d{x},
\qquad \text{where} \qquad
Q = 
\begin{bmatrix}
\frac{1}{\rho(x)} & 0 \\
0 & EI(x) \\
\end{bmatrix}.
\end{equation}

The co-energy variables are found by computing the variational derivative of the Hamiltonian:
\begin{equation}
\begin{aligned}
e_{w} &:= \diffd{H}{\alpha_w} = \diffp{w}{t}(x,t) ,  &\quad \text{Vertical velocity}, \\
e_{\kappa} &:= \diffd{H}{\alpha_{\kappa}} =EI(x) \diffp[2]{w}{x}(x,t),  &\quad \text{Flexural momentum}. \\
\end{aligned}
\end{equation}
Those variables are again collected in vector $\bm{e} = (e_{w}, \, e_{\kappa})^T $, so that the underlying interconnection structure is then found to be:
\begin{equation}
\label{eq:PH_Timo}
\diffp{\bm{\alpha}}{t} = J \mathbf{e},  	\qquad \text{where} \qquad
J = 
\begin{bmatrix}
0 & -\diffp[2]{}{x} \\
\diffp[2]{}{x} & 0 \\
\end{bmatrix}.
\end{equation}

For an infinite-dimensional system, boundary variables have to be defined as well. Those can be found by evaluating the energy rate flow across the boundary. One possible choice among others (see \cite{articleFlavio} for a more exhaustive explanation) for this model is the following:
\begin{equation}
\label{eq:BoundPortEB}
\bm{f}_{\partial} = 
\begin{pmatrix}
e_w(0) \vspace{3pt} \\
\displaystyle \diffp{e_{w}}{x}(0) \vspace{3pt}\\
\displaystyle \diffp{e_{\kappa}}{x}(L) \vspace{3pt}\\
e_{\kappa}(L) \\
\end{pmatrix}, \qquad
\bm{e}_{\partial} = 
\begin{pmatrix}
\displaystyle \diffp{e_{\kappa}}{x}(0) \vspace{1pt}\\
-e_{\kappa}(0) \\
-e_w(L) \\
\displaystyle \diffp{e_{w}}{x}(L)\\
\end{pmatrix}.
\end{equation}
The power flow is then easily evaluated as:
\begin{equation}
\frac{d}{dt} H(t) = \int_{0}^{L} \diffp{\bm{\alpha}}{t} \cdot \bm{e} \, \d{x} = \left\langle \bm{e}_{\partial}, \bm{f}_{\partial} \right\rangle_{{\rm I\!R}^4}.
\end{equation}

The flow variables can now be defined as $\bm{f} = - \diffp{\bm{\alpha}}{t}$, so that the flow space is given by the tuples
$(\bm{f},\bm{f}_{\partial}) \in \mathcal{F}$. Equivalently the effort space is given by $(\bm{e},\bm{e}_{\partial}) \in \mathcal{E}$. The bond space is therefore the Cartesian product of these two spaces:
\begin{equation}
\label{eq:BondEB}
\mathcal{B} := \left\{(\bm{f},\bm{f}_{\partial},\bm{e},\bm{e}_{\partial}) \in \mathcal{F} \times \mathcal{E} \right\}.
\end{equation}

The duality pairing between elements of $\mathcal{B}$ is then defined as follows:
\begin{equation}
\label{eq:bil_Timo}
\ll ((\bm{f}^a,\bm{f}_{\partial}^a), (\bm{e}^a,\bm{e}_{\partial}^a)), ((\bm{f}^b,\bm{f}_{\partial}^b), (\bm{e}^b,\bm{e}_{\partial}^b)) \gg := \int_{0}^{L}\left\{ (\bm{f}^{a})^T \bm{e}^b + (\bm{f}^{b})^T \bm{e}^a \right\} \d{x} + (\bm{f}_{\partial}^{a})^T \bm{e}_{\partial}^b + (\bm{f}_{\partial}^{b}) ^T \bm{f}_{\partial}^a.
\end{equation}
The Stokes-Dirac structure for the Euler-Bernoulli beam is therefore:
\begin{theorem}[From \cite{LeGorrec2005}, Stokes-Dirac structure for the Bernoulli beam]
	Consider the space of power variables {$\mathcal{B}$ defined in \eqref{eq:BondEB}} and the bilinear form (+pairing operator) $\ll , \gg$ given by \eqref{eq:bil_Timo}. Define the following linear subspace $\mathcal{D} \subset \mathcal{F} \times \mathcal{E}$:
	\begin{equation}
	\mathcal{D} =  \left\{(\bm{f},\bm{f}_{\partial},\bm{e},\bm{e}_{\partial}) \in \mathcal{F} \times \mathcal{E} | \; \bm{f} = -J\bm{e} \right\},
	\end{equation}
	{where $\bm{f}_{\partial}$ and $\bm{e}_{\partial}$} were defined in \eqref{eq:BoundPortEB}. Then, it holds $\mathcal{D} = \mathcal{D}^{\perp}$, where $\mathcal{D}^{\perp}$ is understood in the sense of orthogonality with respect to the bilinear product $\ll , \gg$,  i.e $\mathcal{D}$ is a Stokes-Dirac structure. 
\end{theorem}
{
	\begin{remark}
		For what concerns the use of this model for control and simulation purposes, the reader can refer to \cite{stabBeam} for a stability and stabilization proof of the Euler-Bernoulli beam or to \cite{aoues:hal-01738092} for an illustration of a rotating spacecraft with flexible appendages model as PH Bernoulli beams.
	\end{remark}
}
\section{Kirchhoff-Love theory for thin plates}
In this section the classical variational approach (Hamilton's principle) to derive the equation of motions is first detailed. The physical quantities involved and the different energies, of utmost importance for the PH formalism, are reminded.

\subsection{Model and associated variational formulation}
\label{sec:VarKir}
\begin{figure}[tb]
	\centering
	\includegraphics[width=0.8\textwidth]{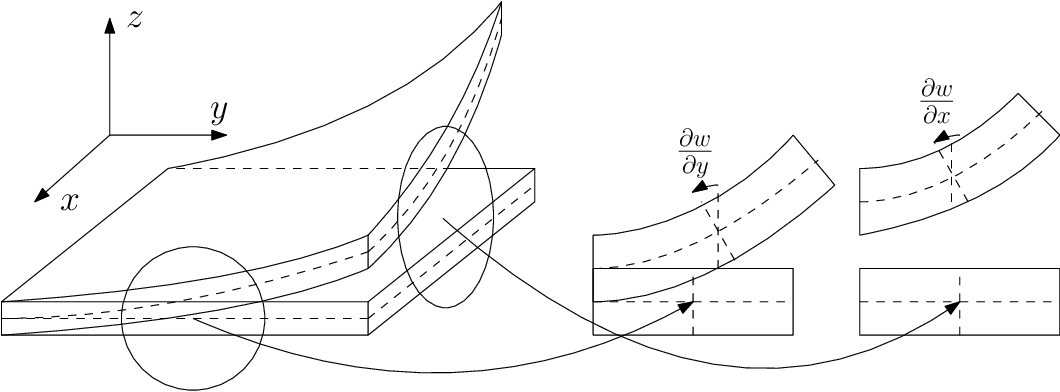}
	\caption{Kinematic assumption for the Kirchhoff plate}
	\label{fig:Kirchh_sketch}
\end{figure}

The Kirchhoff-Love plate formulation rests on the hypothesis of small thickness compared to the in plane dimensions. The notations and symbols are borrowed form \cite{FEM_Cook} and \cite{Oslo}. The displacement field and the strains are defined by assuming that fibers orthogonal to the middle plane remain orthogonal (see Fig. \ref{fig:Kirchh_sketch}). This leads to the following relations for the displacement field
\begin{equation}
u(x,y,z) = -z \diffp{w}{x}, \qquad v(x,y,z) = -z \diffp{w}{y}, \qquad 
w(x,y,z) = w(x,y)
\end{equation}
and for the strains
\begin{equation}
\bm{\epsilon} =  
\begin{pmatrix}
\epsilon_{xx} \\
\epsilon_{yy} \\
\gamma_{xy} \\
\end{pmatrix}  = 
\begin{pmatrix}
\diffp{}{x} & 0 \\
0 & \diffp{}{y} \\
\diffp{}{y} & \diffp{}{x} \\
\end{pmatrix}
\begin{pmatrix}
-z \diffp{w}{x}\\
-z \diffp{w}{y}\\
\end{pmatrix} = 
-z
\begin{pmatrix}
\diffp[2]{w}{x} \\
\diffp[2]{w}{y} \\
2 \diffp{w}{x,y} \\
\end{pmatrix}. 
\end{equation}
The curvature vector is defined as:
\begin{equation}
\label{eq:CurVec}
\bm{\kappa} = 
\begin{pmatrix}
\kappa_{xx} \\
\kappa_{yy} \\
\kappa_{xy} \\
\end{pmatrix} = 
\begin{pmatrix}
\diffp[2]{w}{x} \\
\diffp[2]{w}{y} \\
2 \diffp{w}{x,y} \\
\end{pmatrix}.
\end{equation}

Hooke's constitutive law for isotropic material is considered for the constitutive relation:
\begin{equation}
\bm{\sigma} = 
\bm{E} \bm{\epsilon} , \qquad \bm{E} :=
\frac{E}{1 - \nu^2}
\begin{bmatrix}
1 & \nu & 0 \\
\nu & 1 & 0 \\
0  &  0 & \frac{1-\nu}{2}\\
\end{bmatrix}.
\end{equation}
where $\nu$ is Poisson's ratio and $E$ Young's modulus. {These physical parameters may be inhomogeneous, i.e. $\nu = \nu(x, y, z), \, E = E(x, y, z)$. The generalized momenta are found by integrating the stresses along the fiber:
	\begin{equation*}
	\bm{M} =
	\begin{pmatrix}
	M_{xx} \\
	M_{yy} \\
	M_{xy} \\
	\end{pmatrix} = \left(\int_{-\frac{h}{2} }^{\frac{h}{2}} \bm{E} z^2 \d{z} \right) \; \bm{\kappa},
	\end{equation*}
	where $h$ is the plate thickness. The relation between momenta and curvatures is expressed by the bending rigidity matrix $\bm{D}$:
	\begin{equation}
	\label{eq:BenRigVec}
	\bm{M} = \bm{D} \bm{\kappa} \qquad  \qquad  \bm{D} := \int_{-\frac{h}{2} }^{\frac{h}{2}} \bm{E} z^2 \d{z}.
	\end{equation}
	Now the classical Kirchhoff-Love model for thin plates can be recalled \cite{timoshenko1959theory}:
	\begin{equation}
	\label{eq:KirchhClass}
	\mu \diffp[2]{w}{t}  + \diffp[2]{M_{xx}}{x} + 2 \diffp{M_{xy}}{x,y} + \diffp[2]{M_{yy}}{y}= 0,
	\end{equation}
	where $\mu = \rho h$ is the surface density and $\rho$ the mass density. If the $E$ and $\nu$ coefficients are constant, then the ruling PDE becomes:
	\begin{equation}
	\mu \diffp[2]{w}{t}  + D \Delta^2 w= p ,
	\end{equation}
	where $\Delta^2 = \diffp[4]{}{x} + 2 \diffp[2,2]{}{x,y} + \diffp[4]{}{y}$ is the biLaplacian and $D = \frac{E h^3}{12 (1 - \nu^2)}$ is the bending rigidity modulus.
}
The kinetic and potential energy densities per unit area $\mathcal{K}$ and $\mathcal{U}$, are respectively given by:
\begin{equation*}
\mathcal{K} = \frac{1}{2} \mu \left(\diffp{w}{t}\right)^2, \qquad 
\mathcal{U} = \frac{1}{2} \bm{M} \cdot \bm{\kappa} \, .	
\end{equation*} 

The total energy density is split into kinetic and potential energy
\begin{equation}
\mathcal{H} = \mathcal{K} + \mathcal{U},
\end{equation}
and the corresponding total energies given by the following relations:
\begin{equation}
H = \int_{\Omega} \mathcal{H} \ \d{\Omega}, \qquad K = \int_{\Omega} \mathcal{K} \ \d{\Omega}, \qquad U = \int_{\Omega} \mathcal{U} \ \d{\Omega}.
\end{equation}

\section{PH formulation of the Kirchhoff plate}
In this section the port-Hamiltonian formulation of the Kirchhoff plate is presented first in vectorial form  in \ref{sec:PH_vec_Kir} and then in tensorial form in \ref{sec:PH_ten_Kir}.

\subsection{PH vectorial formulation of the Kirchhoff plate}
\label{sec:PH_vec_Kir}

To obtain a port-Hamiltonian system (PHs) the energy variables as well as the underlying Stokes-Dirac structure, associated with the skew-adjoint operator $J$, have to be properly defined. {Consider the Hamiltonian energy:
	\begin{equation}
	\begin{aligned}
	H &=  \int_{\Omega} \frac{1}{2} \left\{ \mu \left(\diffp{w}{t}\right)^2 + \bm{M} \cdot \bm{\kappa} \right\} \d{\Omega} \\
	&=  \int_{\Omega} \frac{1}{2} \left\{ \mu \left(\diffp{w}{t}\right)^2 + \bm{\kappa}^T \bm{D} \bm{\kappa} \right\} \d{\Omega}.
	\end{aligned}
	\end{equation} 
	The energy variables are then selected to be the linear momentum $\mu \diffp{w}{t}$ and the curvatures $\bm\kappa$, in an analogous fashion with respect to the one-dimensional counterpart of this model, the Euler-Bernoulli beam.} The energy variables are collected in  vector
\begin{equation}
\bm{\alpha} := (\mu w_t,\ \kappa_{xx},\ \kappa_{yy},\ \kappa_{xy})^T,
\end{equation}
where $w_t = \diffp{w}{t}$. The Hamiltonian density is given by the following expression:
\begin{equation}
\label{eq:H_kir}
\mathcal{H} = \frac{1}{2} \bm{\alpha}^T \begin{bmatrix}
\frac{1}{\mu} & 0 \\
0 & \bm{D} \\
\end{bmatrix} \bm{\alpha},  \qquad H = \int_{\Omega} \mathcal{H} \ \d{\Omega}.
\end{equation}
So its variational derivative provides as co-energy variables:
\begin{equation}
\mathbf{e} := \diffd{H}{\bm{\alpha}} = (w_t,\ M_{xx},\ M_{yy},\ M_{xy})^T,
\end{equation}
The port-Hamiltonian system and skew-symmetric operator relating energy and co-energy variables are found to be:
\begin{equation}
\label{eq:PH_Kirchh}
\diffp{\bm{\alpha}}{t} = J \mathbf{e} \quad \text{and} \quad J := 
\begin{bmatrix}
0 & -\diffp[2]{}{x} & -\diffp[2]{}{y} & - \left(\diffp{}{x,y} + \diffp{}{y,x} \right)\\
\diffp[2]{}{x} & 0 & 0 & 0 \\
\diffp[2]{}{y} & 0 & 0 & 0 \\
\diffp{}{x,y} +  \diffp{}{y,x} & 0 & 0 & 0 \\
\end{bmatrix}.
\end{equation}
The first line of the skew-symmetric operator in \eqref{eq:PH_Kirchh} is found by considering Eq. \eqref{eq:KirchhClass}. The remaining lines express Clairaut’s theorem for the vertical displacement. This theorem states that, for smooth functions, higher order partial derivative commute. 
\begin{remark}
	From the Schwarz theorem for $C^2$ functions the mixed derivative could be be expressed as $2 \diffp{}{x,y}$, instead of $\diffp{}{y,x} + \diffp{}{x,y}$. However, in this way the symmetry intrinsically present in $\kappa_{xy} = \diffp{w}{y,x} + \diffp{w}{x,y}$ would be lost. The mixed derivative is here split to reestablish the symmetric nature of curvatures and momenta (that are of tensorial nature as explained in Section \ref{sec:PH_ten_Kir}).
\end{remark}
The boundary variables are obtained by evaluating the time derivative of the Hamiltonian:
\begin{align*}
\dot{H} &= \int_\Omega \diffd{H}{\bm{\alpha}}   \cdot \diffp{\bm{\alpha}}{t} \; \d{\Omega} \\
&= \int_\Omega \left\{ w_t \left(-\diffp[2]{M_{xx}}{x} -\diffp[2]{M_{xx}}{y} - 2 \diffp{M_{xy}}{x,y} \right) + M_{xx} \diffp[2]{w_t}{x} + M_{yy} \diffp[2]{w_t}{y} + 2 M_{xy} \diffp{w_t}{x,y} \right\} \d{\Omega} \\
&= \int_\Omega \left\{ e_1 \left(-\diffp[2]{e_2}{x} -\diffp[2]{e_3}{y} - 2 \diffp{e_4}{x,y}\right) + e_2 \diffp[2]{e_1}{x} + e_3 \diffp[2]{e_1}{y} + 2 e_4 \diffp{e_1}{x,y}  \right\} \d{\Omega}
\end{align*}

\begin{figure}[t]
	\centering
	\includegraphics[height=0.2\textheight]{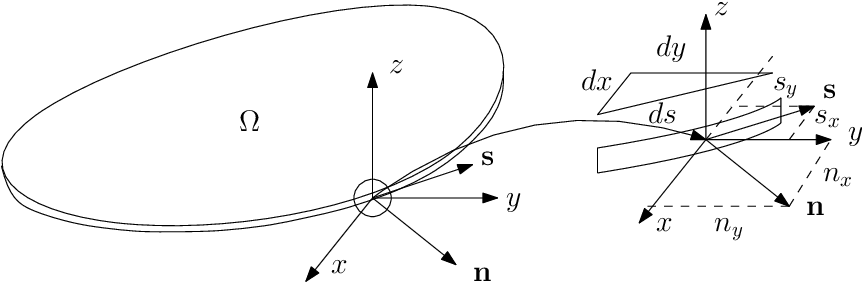}
	\caption{Reference frames and notations.}
	\label{fig:plate_ref}
\end{figure}

In Fig. \ref{fig:plate_ref} the notations for the different reference frames are introduced. By applying Green theorem, considering the split mixed derivative ($2 \diffp{}{x,y} = \diffp{}{x,y} +  \diffp{}{y,x}$):
\begin{align}
\dot{H} = \int_{\partial \Omega}  & \left\{  n_x \left(e_2 \diffp{e_1}{x}  + e_4 \diffp{e_1}{y}  - e_1 \diffp{e_2}{x} - e_1 \diffp{e_4}{y}\right)
\right. \nonumber \\
&   \left. +n_y \left(e_3 \diffp{e_1}{y} + e_4 \diffp{e_1}{x} - e_1 \diffp{e_3}{y} - e_1 \diffp{e_4}{x} \right) \right\} \d{s}.
\end{align}
where $n_x, n_y$ are the components along the $x-$ and the $y-$axis of the normal to the boundary. The variable of integration $s$ is now the curvilinear abscissa which runs along the boundary.

If the physical variables are introduced, then  
\begin{align}
\dot{H} = \int_{\partial \Omega}  & \left\{  n_x \left(M_{xx} \diffp{w_t}{x} + M_{xy} \diffp{w_t}{y} - w_t \, \diffp{M_{xx} }{x}   - w_t \, \diffp{M_{xy}}{y}\right)
\right.  \nonumber\\
&  \left. +n_y \left(M_{yy} \diffp{w_t}{y} + M_{xy} \diffp{w_t}{x} - w_t \, \diffp{M_{yy}}{y} - w_t\, \diffp{M_{xy}}{x} \right) \right\} \d{s}.
\end{align}

\begin{figure}
	\centering
	\includegraphics[width=0.8\textwidth]{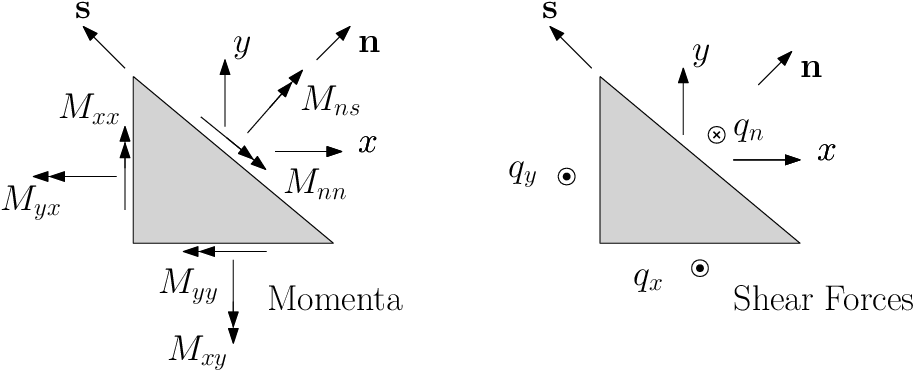}
	\caption{Cauchy law for momenta and forces at the boundary.}
	\label{fig:Cauchy_law}
\end{figure}

Now the following quantities, represented in Fig. \ref{fig:Cauchy_law}, are defined: 
\begin{equation}
\label{eq:QnMnnMns}
\begin{aligned}
\text{Shear Force}& \; \; \quad &q_{n} &:= n_x q_x + n_y q_y,   \\
\text{Flexural momentum}& \quad 
&M_{nn} &:= \bm{n}^T	
\begin{pmatrix}
M_{xx} n_x + M_{xy} n_y \\
M_{xy} n_x + M_{yy} n_y \\
\end{pmatrix}, \\
\text{Torsional momentum}& \quad &M_{ns} &:= \bm{s}^T	
\begin{pmatrix}
M_{xx} n_x + M_{xy} n_y \\
M_{xy} n_x + M_{yy} n_y \\
\end{pmatrix}, 
\end{aligned} \qquad
\begin{aligned}
\bm{n} &= 
\begin{pmatrix}
n_x \\
n_y \\
\end{pmatrix}, \\
\bm{s} &= 
\begin{pmatrix}
-n_y \\
n_x \\
\end{pmatrix},
\end{aligned}
\end{equation}
where $q_x = - \diffp{M_{xx}}{x} -\diffp{M_{xy}}{y}$ and $\; q_y = - \diffp{M_{yy}}{y} -\diffp{M_{xy}}{x}$. The gradient of the vertical velocity can be projected upon the normal and tangential directions to the boundary:
\begin{equation}
\label{eq:gr_dec_ns}
\nabla w_t = \left(\nabla w_t \cdot \bm{n} \right) \,\bm{n} + \left(\nabla w_t \cdot \bm{s} \right) \, \bm{s} = \diffp{w_t}{n} \; \bm{n} +   \diffp{w_t}{s} \; \bm{s}.
\end{equation}
So the time derivative of the Hamiltonian can be finally written as:
\begin{equation}
\dot{H} = \int_{\partial \Omega} \left\{ w_t \, q_n + \diffp{w_t}{s} \, M_{ns} + \diffp{w_t}{n} \, M_{nn}\right\} \d{s}.
\end{equation}
{Variables $w_t$ and $\diffp{w_t}{s}$ are not independent as they are differentially related with respect to derivation along $s$, the curvilinear abscissa of the boundary domain (see for instance \cite{timoshenko1959theory})}. Another integration by part is needed to highlight appropriate {independent} power conjugated variables. Let us suppose that the boundary is a closed and regular curve. Then the integration by parts along a closed boundary leads to:
\begin{equation}
\int_{\partial \Omega} \diffp{w_t}{s} \, M_{ns} \ \d{s}= - \int_{\partial \Omega} \diffp{M_{ns}}{s} \, w_t \ \d{s}.
\end{equation}
The energy balance can be finally written as:
\begin{equation}
\label{eq:energyBal_Kir}
\dot{H} = \int_{\partial \Omega} \left\{ w_t \, \widetilde{q}_n + \diffp{w_t}{n} \, M_{nn}\right\} \d{s}
\end{equation} 
where $\widetilde{q}_n := q_n - \diffp{M_{ns}}{s}$ is the effective shear force.
Equation \eqref{eq:energyBal_Kir} is of utmost importance, since it contains the boundary variables that will be present in the Stokes-Dirac structure defining the port-Hamiltonian system.

\subsubsection{Underlying Stokes-Dirac structure} 

Let $\mathcal{F}$ denote the flow space and let $\mathcal{E}$ denote the effort space. For simplicity we take  $\mathcal{F} \equiv  \mathcal{E} = \mathcal{C}^\infty(\Omega, \mathbb{R}^4)$, the space of smooth vector-valued functions in $\mathbb{R}^4$. Equation \eqref{eq:energyBal_Kir} allows identifying the boundary terms of the underlying Stokes-Dirac structures. The space of boundary {variables} is a vector of four components given by:
\begin{equation*}
\mathcal{Z} = \{\bm{z} | \, \bm{z} = B_{\partial}(\bm{e}), \forall \, \bm{e} \in \mathcal{E} \},  \qquad \bm{z} = \left( \widetilde{q}_n, w_t, M_{nn}, \diffp{w_t}{n} \right)^T .
\end{equation*}

In the case where the differential $J$ operator of order one, the $B_{\partial}$ operator is  a linear operator over the trace of the effort variables.. Here, since the differential $J$ is of order two,  $B_{\partial}$ contains the normal and tangential derivatives at the boundary and so more regularity is required for the boundary variables.
\begin{remark}
	This fact was already stated for 1-D systems in \cite{LeGorrec2005}; here it is the extension to  2-D system with a second-order differential operator $J$.
\end{remark}
This operator reads:
\begin{multline}
B_{\partial}(\bm{e}) = 
\begin{bmatrix}
0 & 0 & 0 & 0 \\
1 & 0 & 0 & 0 \\
0 & n_x^2 & n_y^2 & 2 n_x n_y \\
0 & 0 & 0 & 0 \\
\end{bmatrix}\bm{e}  - 
\begin{bmatrix}
0 & n_x  & 0 & n_y  \\
0 & 0 & 0 & 0 \\
0 & 0 & 0 & 0 \\
0 & 0 & 0 & 0 \\
\end{bmatrix}\diffp{\bm{e}}{x} -
\begin{bmatrix}
0 & 0 & n_y & n_x  \\
0 & 0 & 0 & 0 \\
0 & 0 & 0 & 0 \\
0 & 0 & 0 & 0 \\
\end{bmatrix}\diffp{\bm{e}}{y} \\
+ \diffp{}{n}\left(
\begin{bmatrix}
0 & 0 & 0 & 0 \\
0 & 0 & 0 & 0 \\
0 & 0 & 0 & 0 \\
1 & 0 & 0 & 0 \\
\end{bmatrix} \bm{e} \right) - \diffp{}{s}\left(
\begin{bmatrix}
0 & - n_x n_y & n_x n_y & n_x^2-n_y^2 \\
0 & 0 & 0 & 0 \\
0 & 0 & 0 & 0 \\
0 & 0 & 0 & 0 \\
\end{bmatrix} \bm{e} \right).
\end{multline}

\begin{theorem}[Stokes-Dirac structure fr the Kirchhoff Plate]
	The set
	\begin{equation}
	\mathcal{D} := \left\{ (\bm{f},\bm{e},\bm{z}) \in \mathcal{F}\times\mathcal{E}\times\mathcal{Z} \; | \; \bm{f}= - \diffp{\bm{\alpha}}{t} = -J \bm{e}, \; \bm{z} = B_{\partial}(\bm{e}) \right\}
	\end{equation} 
	is a Stokes-Dirac structure with respect to the pairing
	\begin{equation}
	\ll (\bm{f}_1, \bm{e}_1, \bm{z}_1), (\bm{f}_2, \bm{e}_2, \bm{z}_2) \gg  \,= \int_{\Omega} \left[ \bm{e}_1^T \bm{f}_2 + \bm{e}_2^T \bm{f}_1 \right] \d{\Omega}  + \int_{\partial \Omega} B_J(\bm{z}_1, \bm{z}_2) \, \d{s},
	\end{equation}
	where $B_J$ is a symmetric operator, arising from a double application of the Green theorem. It reads
	\begin{equation}
	\begin{aligned}
	B_J(\bm{z}_1,\bm{z}_2) = \, &\widetilde{q}_{n, 2} \ w_{t, 1} + M_{nn, 2} \ \diffp{w_{t, 1}}{n} \\
	+ \, &\widetilde{q}_{n, 1} \ w_{t, 2} + M_{nn, 1} \ \diffp{w_{t, 2}}{n} \\
	= &\bm{z}_1^T \, B_J \, \bm{z}_2 \\
	\end{aligned} \, , \qquad	B_J = 
	\begin{bmatrix}
	0 & 1 & 0 & 0 \\
	1 & 0 & 0 & 0 \\
	0 & 0 & 0 & 1 \\
	0 & 0 & 1 & 0 \\ 
	\end{bmatrix}.
	\end{equation}
\end{theorem}
\begin{proof}
	A regular boundary will be assumed for the proof.
	\paragraph{\textbf{Step I}} 
	The first implication is $\mathbb{D} \subseteq \mathbb{D}^T$. This is true if $ \forall \;\bm\omega_\alpha = (\bm{f}_\alpha, \bm{e}_\alpha, \bm{z}_\alpha)$ and $\bm\omega_\beta = (\bm{f}_\beta, \bm{e}_\beta, \bm{z}_\beta) \in \mathbb{D}$ then $\ll \bm\omega_\alpha, \bm\omega_\beta \gg = 0$. The integral over the domain reads:
	\begin{align*}
	&\int_{\Omega}  \left[ \bm{e}_\alpha^T \bm{f}_\beta + \bm{e}_\beta^T \bm{f}_\alpha \right] d\Omega
	= \int_{\Omega} \left\{ e_1^\alpha \left( \diffp[2]{e_2^\beta}{x} + \diffp[2]{e_3^\beta}{y} + \diffp{e_4^\beta}{x,y} \right) - e_2^\alpha \diffp[2]{e_1^\beta}{x} -e_3^\alpha \diffp[2]{e_1^\beta}{y} \right.\\ 
	&  - 2 e_4^\alpha \diffp{e_1^\beta}{x,y} + 
	\left. e_1^\beta \left( \diffp[2]{e_2^\alpha}{x} + \diffp[2]{e_3^\alpha}{y} + \diffp{e_4^\alpha}{x,y} \right) - e_2^\beta \diffp[2]{e_1^\alpha}{x} -e_3^\beta \diffp[2]{e_1^\alpha}{y} - 2 e_4^\beta \diffp{e_1^\alpha}{x,y} \right\} \d{\Omega}. \\
	\end{align*}
	Once the Green theorem has been applied,  the relevant quantities expressed by equations~\eqref{eq:QnMnnMns} pop up as follows:
	\begin{multline}
	\int_{\Omega}  \left[ \bm{e}_\alpha^T \bm{f}_\beta + \bm{e}_\beta^T \bm{f}_\alpha \right] d\Omega = \int_{\partial \Omega}  \left\{ e_1^\alpha \left( \diffp{e_2^\beta}{x} n_x + \diffp{e_3^\beta}{y} n_y +\diffp{e_4^\beta}{y} n_x + \diffp{e_4^\beta}{x} n_y \right) \right. \\
	\left.   + e_1^\beta \left( \diffp{e_2^\alpha}{x} n_x + \diffp{e_3^\alpha}{y} n_y +\diffp{e_4^\alpha}{y} n_x + \diffp{e_4^\alpha}{x} n_y \right) - e_2^\alpha \diffp{e_1^\beta}{x} n_x - e_2^\beta \diffp{e_1^\alpha}{x} n_x   \right.\\
	\left. - e_3^\alpha \diffp{e_1^\beta}{y} n_y  - e_3^\beta \diffp{e_1^\alpha}{y} n_y  - e_4^\alpha \left( \diffp{e_1^\beta}{y} n_x + \diffp{e_1^\beta}{x} n_y \right) - e_4^\beta \left( \diffp{e_1^\alpha}{y} n_x + \diffp{e_1^\alpha}{x} n_y \right) \right\} \d{s}. 
	\end{multline}
	
	Moreover only the kinematically independent quantities have to be considered, leading to the final result:
	\begin{equation}
	- \int_{\partial \Omega} \left( w_t^\alpha \widetilde{q}_n^\beta + w_t^\beta \widetilde{q}_n^\alpha + \diffp{w_t^\alpha}{n} \, M_{nn}^\beta + \diffp{w_t^\beta}{n} \, M_{nn}^\alpha \right) ds =  - \int_{\partial \Omega}  B_J(\bm{z}_\alpha, \bm{z}_\beta) ds.
	\end{equation}		
	This concludes the first part of the proof.
	\paragraph{\textbf{Step II}}
	For the second implication, i.e. $\mathcal{D}^{\perp} \subseteq \mathcal{D}$. Let us take $\bm\omega_\alpha \in \mathcal{D}^{\perp}, \forall \, \bm\omega_\beta \in \mathcal{D}$. Then the bilinear form, once the Green theorem has been applied, provides the following
	\begin{multline}
	\int_{\Omega} \left\{ e_1^\beta \left( f_1^\alpha - \diffp[2]{e_2^\alpha}{x} - \diffp[2]{e_3^\alpha}{y} - 2 \diffp{e_4^\alpha}{x,y}  \right) + e_2^\beta \left( f_2^\alpha + \diffp[2]{e_1^\alpha}{x}  \right) + e_3^\beta \left( f_3^\alpha + \diffp[2]{e_1^\alpha}{y}  \right) + \right. \\
	\left. e_4^\beta \left( f_4^\alpha + 2 \diffp{e_1^\alpha}{x,y}  \right)  \right\} d\Omega + \int_{\partial \Omega} \left\{ w_t^\beta \left( \diffp{e_2^\alpha}{x} n_x + \diffp{e_3^\alpha}{y} n_y + \diffp{e_4^\alpha}{y} n_x + \diffp{e_4^\alpha}{y} n_y  \right)  \right. \\
	\left.  - q_n^\beta e_1^\alpha - \diffp{w_t^\beta}{x} e_2^\alpha n_x - \diffp{w_t^\beta}{y} e_3^\alpha n_y   - e_4^\alpha \left( \diffp{w_t^\beta}{y} n_x + \diffp{w_t^\beta}{x} n_y  \right) - \diffp{e_1^\alpha}{x} M_{xx}^\beta n_x  - \diffp{e_1^\alpha}{y} M_{yy}^\beta n_y \right.  \\ 
	\left. - M_{xy}^\beta \left(\diffp{e_1^\alpha}{y} n_x + \diffp{e_1^\alpha}{x} n_y  \right) + \widetilde{q}_n^\beta z_1^\alpha + w_t^\beta z_2^\alpha + M_{nn}^\beta z_3^\alpha + \diffp{w_t^\beta}{n} z_4^\beta \right\} ds = 0. 
	\end{multline}
	Since the relation has to be valid for each $\bm\omega_\beta \in \mathcal{D}$ the flux and effort variables are in $\mathcal{D}$. For the boundary terms the same procedure as before has to be applied by considering the definition of the momenta over the boundary (see equation \eqref{eq:QnMnnMns}). Then it can be stated that $\bm\omega_\alpha \in \mathcal{D}$.
\end{proof}
\subsubsection{Including dissipation and external forces in the model}
Distributed forces or control and dissipative relations can be easily included in an augmented Stokes-Dirac structure by simply defining the appropriate conjugated variables. \\

If distributed forces have to be considered, then the set
\begin{multline}
\mathcal{D}_d := \Bigl\{ (\bm{f}, \bm{f}_d, \bm{e}, \bm{e}_d, \bm{z}) \in \mathcal{F}\times\mathcal{F}_d\times\mathcal{E}\times\mathcal{E}_d\times\mathcal{Z} \; |  \\ \bm{f}= - \diffp{\bm{\alpha}}{t} = -J \bm{e} - G_d \bm{f}_d, \; \bm{e}_d = G^*_d \bm{e}, \;  \bm{z} = B_{\partial}(\bm{e}) \Bigr\}
\end{multline}
is a Stokes-Dirac structure with respect to the paring 
\begin{multline}
\ll (\bm{f}_1, \bm{f}_{d, 1} \bm{e}_1, \bm{e}_{d, 1}, \bm{z}_1), (\bm{f}_2, \bm{f}_{d, 2} \bm{e}_2, \bm{e}_{d, 2}, \bm{z}_2) \gg  \,= \\
\int_{\Omega} \left[ \bm{e}_1^T \bm{f}_2 + \bm{e}_2^T \bm{f}_1 + \bm{e}_{d, 1}^T \bm{f}_{d, 2} + \bm{e}_{d, 2}^T \bm{f}_{d, 1} \right] \, d\Omega  + \int_{\partial \Omega} B_J(\bm{z}_1, \bm{z}_2) \, ds.
\end{multline}
If gravity has to be included, then $G_d=[1, 0, 0, 0]^T, f_d = -\mu g$.  \\

Analogously dissipation can be included in an augmented Dirac structure. As an example, the ruling PDE, once a dissipative term of fluid damping type is considered, reads:
\begin{equation}
\mu \diffp[2]{w}{t}  + r \diffp{w}{t}+ \diffp[2]{M_{xx}}{x} + 2 \diffp{M_{xy}}{x,y} + \diffp[2]{M_{yy}}{y} = 0,
\end{equation}
where $r>0$ is the damping coefficient.
If this equation is rewritten using the port-Hamiltonian formalism then we get:
\begin{equation}
\diffp{\bm\alpha}{t} = \left(J - R \right) \bm{e}, \qquad 
R := 
\begin{bmatrix}
r & 0 & 0 & 0 \\
0 & 0 & 0 & 0 \\
0 & 0 & 0 & 0 \\
0 & 0 & 0 & 0 \\
\end{bmatrix}.
\end{equation}
The $R$ matrix, which is a symmetric, semi-positive definite operator, can be decomposed as
\begin{equation}
R = G_R S G_R^* , 
\end{equation}
where $S=r$ is a coercive operator (in this case simply a positive scalar), $G_R =  \left(1 \ 0 \ 0 \ 0 \right)^T$ and $G_R^*$ denotes the adjoint operator to $G_R$. The augmented structure
\begin{multline}
\mathcal{D}_r := \Bigl\{ (\bm{f}, \bm{f}_r) \in \mathcal{F}, \ (\bm{e}, \bm{e}_r) \in \mathcal{E}, \ \bm{z} \in \mathcal{Z} \; |  \\
\bm{f} = - \diffp{\bm{\alpha}}{t} = -J \bm{e} - G_R \bm{f}_r, \; \bm{f}_r = - S \bm{e}_r, \; \bm{e}_r = G^*_R \bm{e},  \;  \bm{z} = B_{\partial}(\bm{e})	
\Big\}
\end{multline}
is a Stokes-Dirac structure with respect to the paring 
\begin{multline}
\ll (\bm{f}_1, \bm{f}_{r, 1}, \bm{e}_1, \bm{e}_{r, 1}, \bm{z}_1), (\bm{f}_2, \bm{f}_{r, 2}, \bm{e}_2, \bm{e}_{r, 2}, \bm{z}_2) \gg  \,= \\
\int_{\Omega} \left[ \bm{e}_1^T \bm{f}_2 + \bm{e}_2^T \bm{f}_1 + \bm{e}_{r, 1}^T \bm{f}_{r, 2} + \bm{e}_{r, 2}^T \bm{f}_{r, 1} \right] d\Omega  + \int_{\partial \Omega} B_J(\bm{z}_1, \bm{z}_2) \, ds.
\end{multline}

\begin{remark}
	More involved dissipation models can be found in \cite{DissDenis}. More specifically, for Kirchhoff plate, some specific damping models can be found in \cite{LambourgJASA}.
\end{remark}

\subsection{PH tensorial formulation of the Kirchhoff plate}
\label{sec:PH_ten_Kir}

In section \ref{sec:PH_vec_Kir} the Stokes-Dirac structure of the Kirchhoff plate was found by using a vectorial notation for the curvatures and momenta. {In fact} these variables are of tensorial nature and in the following the tensorial formulation takes the place of the vectorial one. First let us rewrite the momenta and curvatures as symmetric matrices (corresponding to the choice of a Cartesian frame for the representation of tensors):
\begin{equation}
\mathbb{K} = 
\begin{bmatrix}
\kappa_{xx} &  \kappa_{xy}\\
\kappa_{xy} & \kappa_{yy} \\
\end{bmatrix}, \qquad
\mathbb{M} =
\begin{bmatrix}
M_{xx} & M_{xy} \\
M_{xy} & M_{yy} \\
\end{bmatrix},
\end{equation}
where now, with a slight abuse of notation, $\kappa_{xy}$ {differs by $1/2$ from the definition given in equation \eqref{eq:CurVec}}, i.e. $\kappa_{xy} =\diffp{w}{x,y}$. All the other quantities stay the same with respect to what stated in section \ref{sec:VarKir}. The Hamiltonian energy is written as:
\begin{equation}
\label{eq:H_Kir_Ten}
H = \int_{\Omega} \left\{ \frac{1}{2} \mu \left(\diffp{w}{t} \right)^2 + \frac{1}{2} \mathbb{M} \cddot \mathbb{K}  \right\}  \d{\Omega} ,
\end{equation}
where the tensor contraction in Cartesian coordinates is expressed as
\[\mathbb{M} \cddot \mathbb{K} = \sum_{i,j = 1}^{2} M_{ij} \kappa_{ij} = \Tr(\mathbb{M}^T \mathbb{K}). \]
For what concerns the choice of the energy variables, a scalar and a tensor variable are considered:
\begin{equation}
\alpha_w = \mu \diffp{w}{t}, \qquad \mathbb{A}_{\kappa} = \mathbb{K}.	\end{equation}
The co-energy variables are found by computing the variational derivative of the Hamiltonian:
\begin{equation}
e_w := \diffd{H}{\alpha_w} = \diffp{w}{t} := w_t,  \qquad  \mathbb{E}_{\kappa} := \diffd{H}{\mathbb{A}_{\kappa}} = \mathbb{M}.
\end{equation}
\begin{remark}
	For the variational derivative with respect to a tensor, see Propostion 1 in \cite{BrugnoliMin}.
\end{remark}
The port-Hamiltonian system {\eqref{eq:PH_Kirchh} is now rewritten as:}
\begin{equation}
\label{eq:PH_sys_Kir_Ten_pr}
\begin{cases}
\displaystyle\diffp{\alpha_w}{t} &= - \mathrm{div}(\mathrm{Div}(\mathbb{E}_{\kappa})), \vspace{1mm} \\
\displaystyle\diffp{\mathbb{A}_{\kappa}}{t} &= \mathrm{Grad}(\mathrm{grad}(e_w)),
\end{cases}
\end{equation}
where $\mathrm{div}$ and $\mathrm{Div}$ denote the divergence of a vector and of a tensor respectively. The operator $\mathrm{Grad}$ denotes the symmetric gradient:
\begin{equation}
\mathrm{Grad}(\bm{a}) =  \frac{1}{2} \left(\nabla \otimes \bm{a} + \left(\nabla \otimes \bm{a}\right)^T \right).
\end{equation}

The operator $\mathrm{Grad} \circ \mathrm{grad}$ corresponds to the Hessian operator. In Cartesian coordinates it reads:
\begin{equation}
\mathrm{Grad} \circ \mathrm{grad} = 
\begin{bmatrix}
\diffp[2]{}{x}  &  \diffp{}{x,y} \\
\diffp{}{y,x}   &  \diffp[2]{}{y} \\
\end{bmatrix}.
\end{equation}

\begin{theorem}
	The operator $\mathrm{Grad} \circ \mathrm{grad}$, corresponding to the Hessian operator, is the adjoint of the double divergence $\mathrm{div} \circ \mathrm{Div}$.
\end{theorem}

\begin{proof}
	Let us consider the Hilbert space of the square integrable symmetric square tensors of size $n \times n$ over an open connected set $\Omega$. This space will be denoted by $\mathscr{H}_1 = L^2(\Omega, \mathbb{R}^{n \times n}_{\text{sym}})$. This space is endowed with the integral of the tensor contraction as scalar product:
	\[\left\langle \mathbb{E} , \mathbb{F} \right\rangle_{\mathscr{H}_1} = \int_{\Omega}  \mathbb{E} \cddot \mathbb{F} \; \d{\Omega} = \int_{\Omega} \Tr(\mathbb{E}^T \mathbb{F}) \; \d{\Omega}, \quad \forall  \mathbb{E} , \mathbb{F} \in [L^2_{\text{sym}}(\Omega)]^{n \times n}. \]
	
	Consider the Hilbert space $\mathscr{H}_2 =L^2(\Omega)$ of scalar square integrable functions, endowed with the inner product:
	\begin{equation}
	\left\langle e,f \right\rangle_{\mathscr{H}_2}  \,= \int_{\Omega} e f \, \d{\Omega}.
	\end{equation} 
	Let us consider the double divergence operator defined as: 
	\[
	\begin{aligned}
	A: \; \mathscr{H}_1& \rightarrow \mathscr{H}_2, \\
	\mathbb{E}& \rightarrow \mathrm{div}(\mathrm{Div}(\mathbb{E})) = e, \\	\end{aligned}
	\qquad \text{with } \bm{e} = \mathrm{div}(\mathrm{Div}(\mathbb{E})) = \sum_{i = 1}^n \sum_{j = 1}^n \diffp{\mathbb{E}_{ij}}{x_i,x_j}.
	\]
	We {shall} identify $A^*$
	\[
	\begin{aligned}
	A^*: \; \mathscr{H}_2& \rightarrow \mathscr{H}_1, \\
	f& \rightarrow  A^* f = \mathbb{F}, \\
	\end{aligned}
	\]
	such that \[
	\left\langle A \mathbb{E} , f \right\rangle_{\mathscr{H}_2} = \left\langle \mathbb{E} , A^* f \right\rangle_{\mathscr{H}_1},
	\begin{aligned} \qquad
	&\forall \,\mathbb{E} \in \mathrm{Domain}(A) \subset \mathscr{H}_1 \\
	&\forall \,f \in \mathrm{Domain}(A^*) \subset \mathscr{H}_2
	\end{aligned}
	\]

	The function have to belong to the operator domain, so for instance $f \in \mathcal{C}_0^2(\Omega) \in \mathrm{Domain}(A^*)$ the space of twice differentiable scalar functions with compact support on an open simply connected set $\Omega$ and additionally $\mathbb{E}$ can be chosen in the set $\mathcal{C}_{0}^2(\Omega, \mathbb{R}^{2 \times 2}_{\text{sym}}) \in \mathrm{Domain}(A)$, the space of twice differentiable $2 \times 2$ symmetric tensors with compact support on $\Omega$. { A classical result is the fact that the adjoint of the vector divergence is $\mathrm{div}^* = -\mathrm{grad}$ as stated in \cite{kurula2012duality}. This may be generalized to the adjoint of the tensor divergence $\mathrm{Div}^* = -\mathrm{Grad}$ (see Theorem 4 of \cite{BrugnoliMin}).} Considering that $A$ is the composition of two different operators $A = \mathrm{div} \circ \mathrm{Div}$ and that the adjoint of a composed operator is the adjoint of each operator in reverse order, i.e. $(B \circ C)^* = C^* \circ B^*$, then it can be stated
	\[
	A^* = (\mathrm{div} \circ \mathrm{Div})^* = \mathrm{Div}^* \circ \mathrm{div}^* = \mathrm{Grad} \circ \mathrm{grad}.
	\]  
	Since only formal adjoints are being looked for, this concludes the proof.
\end{proof}
If the variables in system \eqref{eq:PH_sys_Kir_Ten_pr} are gathered together the formally skew-symmetric operator $J$ can be highlighted:
\begin{equation}
\label{eq:PH_sys_Kir_Ten}
\diffp{}{t}
\begin{pmatrix}
\alpha_w \\
\mathbb{A}_{\kappa} \\
\end{pmatrix} = 
\underbrace{\begin{bmatrix}
	0  &  - \mathrm{div} \circ \mathrm{Div} \\
	\mathrm{Grad} \circ \mathrm{grad} & 0 \\
	\end{bmatrix}}_{J}
\begin{pmatrix}
e_w \\
\mathbb{E}_{\kappa} \\
\end{pmatrix}.
\end{equation}
where all zeros are intended as nullifying operator from the space of input variables to the space of output variables.
\begin{remark}
	The interconnection structure $J$ now resembles that of the Bernoulli beam. The double divergence and the double gradient coincide, in dimension one, with the second derivative.
\end{remark}
Again the boundary {port variables} can be found by evaluating the time derivative of the Hamiltonian:
\begin{equation}
\begin{aligned}
\dot{H}&= \int_{\Omega} \left\{ \diffp{\alpha_w}{t} e_w  + \diffp{\mathbb{A}_{\kappa}}{t} \cddot \mathbb{E}_{\kappa} \right\} \d\Omega \\
&= \int_{\Omega} \left\{ - \mathrm{div}(\mathrm{Div}(\mathbb{E}_{\kappa})) e_w + \mathrm{Grad}(\mathrm{grad}(e_w)) \cddot \mathbb{E}_{\kappa} \right\} \d\Omega, &\qquad \text{Integration by parts} \\
&=  \int_{\partial \Omega} \left\{ \underbrace{  - \bm{n} \cdot \mathrm{Div}(\mathbb{E}_{\kappa}) }_{q_n} e_w + \left[\bm{n} \otimes \mathrm{grad}(e_w) \right]\cddot \mathbb{E}_{\kappa} \right\} \d{s}, &\qquad \text{See \eqref{eq:QnMnnMns} and \eqref{eq:gr_dec_ns} }  \\
&=\int_{\partial \Omega} \left\{ q_n e_w + \diffp{e_w}{n} \underbrace{\left(\bm{n} \otimes \bm{n} \right) \cddot \mathbb{E}_{\kappa}}_{M_{nn}} + \diffp{e_w}{s}  \underbrace{\left(\bm{n} \otimes\bm{s} \right)\cddot \mathbb{E}_{\kappa}}_{M_{ns}}   \right\} \d{s}, &\qquad \text{Dyadic properties}\\
&=\int_{\partial \Omega} \left\{ q_n w_t + \diffp{w_t}{n} M_{nn} + \diffp{w_t}{s}  M_{ns}   \right\} \d{s}. 
\end{aligned}
\end{equation}
{
	\begin{remark}
		The definitions
		\[q_n=  - \bm{n} \cdot \mathrm{Div}(\mathbb{E}_{\kappa}), \quad M_{nn} = \left(\bm{n} \otimes \bm{n} \right) \cddot \mathbb{E}_{\kappa}, \quad M_{ns} = \left(\bm{n} \otimes \bm{s} \right) \cddot \mathbb{E}_{\kappa}\] 
		are exactly the same as those given in \eqref{eq:QnMnnMns}. The tensorial formalism allows a more compact writing.
	\end{remark}
}
The kinematically independent variables must be highlighted. The tangential derivative has to be moved on the torsional momentum. In order to do that, the boundary needs to be split in a collection of regular subsets $\Gamma_{i}$, such that $\partial \Omega = \bigcup_{\Gamma_{i} \subset \partial \Omega} \Gamma_{i}$:
\begin{equation}
\begin{aligned}
\int_{\partial \Omega} \diffp{w_t}{s} \, M_{ns} \ \d{s} &= \sum_{\Gamma_{i} \subset \partial \Omega} \int_{\Gamma_{i}}  \diffp{w_t}{s} \, M_{ns} \ \d{s} \\
&= \sum_{\Gamma_{i} \subset \partial \Omega} \left[ M_{ns} w_t\right]_{\partial \Gamma_i} - \int_{\partial \Omega} \diffp{M_{ns}}{s} \, w_t \ \d{s}.
\end{aligned}
\end{equation}

If a regular boundary is considered the final energy balance is exactly the same as the obtained with the vectorial notation, namely:
\begin{equation}
\dot{H} = \int_{\partial \Omega} \left\{ w_t \, \widetilde{q}_n + \diffp{w_t}{n} \, M_{nn}\right\} \ \d{s},  \qquad \text{where }\widetilde{q}_n := q_n - \diffp{M_{ns}}{s}.
\end{equation} 
The tensorial formulation allows highlighting the intrinsic differential operators. Furthermore the symmetric nature of the variables is explicitly expressed by the usage of symmetric tensors. {Now that the energy balance has been established in terms of the boundary variables the Stokes-Dirac structure for the Kirchhoff plate in tensorial form can be defined. Consider now the bond space:
	\begin{equation}
	\label{eq:bondKir}
	\mathcal{B} := \left\{(\bm{f}, \bm{e}, \bm{z}) \in \mathcal{F} \times \mathcal{E} \times \mathcal{Z} \right\},
	\end{equation}
	where $\mathcal{F}=  \mathscr{L}^2(\Omega) :=  L^2(\Omega) \times L^2(\Omega, \mathbb{R}^{2 \times 2}_{\text{sym}})$ and $ \mathcal{E} =  \mathscr{H}^2(\Omega) = H^{2}(\Omega) \times H^{\text{div Div}}(\Omega, \mathbb{R}^{2 \times 2}_{\text{sym}})$. The space $H^{\text{divDiv}}(\Omega, \mathbb{R}^{2 \times 2}_{\text{sym}})$ is such that
	\[
	H^{\text{div Div}}(\Omega, \mathbb{R}^{2 \times 2}_{\text{sym}}) = \left\{ \mathbb{A} \in L^2(\Omega, \mathbb{R}^{2 \times 2}_{\text{sym}}) \; \vert \; \mathrm{div}(\mathrm{Div}(\mathbb{A})) \in L^2(\Omega) \right\}.
	\] Consider the space of boundary port variables:
	\begin{equation}
	\begin{gathered}
	\mathcal{Z} := \left\{ \bm{z} \; \vert \; \bm{z} = \begin{pmatrix} \bm{f}_{\partial} \\ \bm{e}_{\partial} \end{pmatrix} \right\}, \qquad \text{with} \quad
	\bm{f}_\partial = 
	\begin{pmatrix}
	w_t \\ \diffp{w_t}{n} \\
	\end{pmatrix}, \qquad
	\bm{e}_\partial = 
	\begin{pmatrix}
	\widetilde{q}_n \\ M_{nn} \\
	\end{pmatrix}.
	\end{gathered}
	\end{equation}
	The duality pairing between elements of $\mathcal{B}$ is then defined as follows:
	\begin{equation}
	\label{eq:bil_Kir_Ten}
	\left\langle \left\langle (\bm{f}_1, \bm{e}_1, \bm{z}_1), (\bm{f}_2, \bm{e}_2, \bm{z}_2) \right\rangle \right\rangle :=  \left\langle \bm{e}_1, \bm{f}_2 \right\rangle_{\mathscr{L}^2(\Omega)}  +  \left\langle \bm{e}_2,  \bm{f}_1  \right\rangle_{\mathscr{L}^2(\Omega)} + \int_{\partial \Omega} B_J(\bm{z}_1, \bm{z}_2) \d{s} , 
	\end{equation}
	where the pairing $\left\langle \cdot, \cdot \right\rangle_{\mathscr{L}^2(\Omega)}$ is the $L^2$ inner product on space $\mathscr{L}^2(\Omega)$ and $B_J(\bm{z}_1, \bm{z}_2) := (\bm{f}_{\partial, 1})^T \bm{e}_{\partial, 2} + (\bm{f}_{\partial, 2}) ^T \bm{e}_{\partial, 1}$. 
	\begin{theorem}[Stokes-Dirac Structure for the Kirchhoff plate in tensorial form]
		Consider the space of power variables $\mathcal{B}$ defined in \eqref{eq:bondKir} and the matrix differential operator $J$ in~\eqref{eq:PH_sys_Kir_Ten}. By theorem 2 in \cite{BrugnoliMin} the linear subspace $\mathcal{D} \subset \mathcal{B}$
		\begin{equation}
		\mathcal{D} =  \left\{(\bm{f}, \bm{e},\bm{z}) \in \mathcal{B} | \; \bm{f} = - \diffp{\bm{\alpha}}{t} = -J\bm{e}, \; \bm{z} = \begin{pmatrix} \bm{f}_{\partial} \\ \bm{e}_{\partial} \end{pmatrix} 
		\right\},
		\end{equation}
		is a Stokes-Dirac structure with respect to the pairing $\left\langle \left\langle \cdot, \cdot \right\rangle \right\rangle$ given by \eqref{eq:bil_Kir_Ten}.
	\end{theorem}
}

\section{Discretization of the Kirchhoff plate using a Partioned Finite Element Method} 
Following the procedure illustrated in \cite{CardosoRibeiro2018} the Kirchhoff plate written as a port-Hamiltonian system can be discretized by using a Partitioned Finite Element Method (PFEM). {This method is an extension of the Mixed Finite Element Method to the case of pH systems and requires the integration by parts to be performed, so that the symplectic structure is preserved. It consists of three different steps:
	\begin{enumerate}
		\item the system is first put into weak form;
		\item once the boundary control of interest is selected, the corresponding subsystem is integrated by parts;
		\item the problem is discretized by using a Mixed Finite Element method.  
	\end{enumerate}
	The weak form is illustrated using the tensorial formulation. Two different kind of boundary controls will be shown:
	\begin{enumerate}
		\item boundary control through forces and momenta, in this case the first line of \eqref{eq:PH_sys_Kir_Ten} is integrated by parts (in $\S$ \ref{sec:BC_forces});
		\item boundary control through kinematic variables, in this case the second line of \eqref{eq:PH_sys_Kir_Ten} is integrated by parts (in $\S$ \ref{sec:BC_displ}).
	\end{enumerate}
}

\subsection{Weak form}
The same procedure detailed above can be used on system \eqref{eq:PH_sys_Kir_Ten}. In this case the test functions are of scalar or tensorial nature. Keeping the same notation than in Section \ref{sec:PH_ten_Kir} the scalar test function is denoted by $v_w$, the tensorial one by $\mathbb{V}_{\kappa}$. \\

\subsubsection{Boundary control through forces and momenta}
\label{sec:BC_forces}
The fist line of \eqref{eq:PH_sys_Kir_Ten} is multiplied  by $v_w$ (scalar multiplication), the second line by $\mathbb{V}_{\kappa}$ (tensor contraction).

\begin{align}
\int_{\Omega} v_w \diffp{\alpha_w}{t} \,  \d\Omega &=  \int_{\Omega} -v_w \mathrm{div}(\mathrm{Div}(\mathbb{E}_{\kappa})) \, \d\Omega,  \label{eq:wf1_kir_ten}\\
\int_{\Omega} \mathbb{V}_{\kappa} \cddot \diffp{\mathbb{A}_{\kappa}}{t} \,  \d\Omega &= \int_{\Omega} \mathbb{V}_{\kappa} \cddot  \mathrm{Grad}(\mathrm{grad}(e_w)) \,   \d\Omega.  \label{eq:wf2_kir_ten}
\end{align}

The right hand side of equation \eqref{eq:wf1_kir_ten} has to be integrated by parts twice:
\begin{equation}
\label{eq:line1_wf1_ten}
\int_{\Omega} - v_w \mathrm{div}(\mathrm{Div}(\mathbb{E}_{\kappa})) \, \d\Omega = \int_{\partial \Omega} \underbrace{- \bm{n} \cdot \mathrm{Div}(\mathbb{E}_{\kappa})}_{q_n} v_w \, \d{s} + \int_{\Omega} \mathrm{grad}(v_w)  \cdot \mathrm{Div}(\mathbb{E}_{\kappa}) \, \d\Omega
\end{equation}
Applying again the integration by parts leads to:
\begin{equation}
\label{eq:line1_wf2_ten}
\int_{\Omega} \mathrm{grad}(v_w)  \cdot \mathrm{Div}(\mathbb{E}_{\kappa}) \, \d\Omega = \int_{\partial \Omega} \mathrm{grad}(v_w)  \cdot \left( \bm{n} \cdot \mathbb{E}_{\kappa} \right) \, \d{s} -  \int_{\Omega}\mathrm{Grad}(\mathrm{grad}(v_w))  \cddot \mathbb{E}_{\kappa} \, \d\Omega
\end{equation}
The usual additional manipulation is performed on the boundary term containing the momenta, so that the proper boundary values arise:
\begin{equation}
\label{eq:line1_ipbc_ten}
\begin{aligned}
\int_{\partial \Omega} \mathrm{grad}(v_w)  \cdot \left( \bm{n} \cdot \mathbb{E}_{\kappa} \right) \, \d{s} &= \int_{\partial \Omega} \left(\diffp{v_w}{n} \bm{n} + \diffp{v_w}{s} \bm{s} \right)  \cdot \left( \bm{n} \cdot \mathbb{E}_{\kappa} \right) \, \d{s} \\
&= \int_{\partial \Omega} \left\{  \diffp{v_w}{n}  \underbrace{\left(\bm{n} \otimes \bm{n} \right) \cddot \mathbb{E}_{\kappa} }_{M_{nn}} +  \diffp{v_w}{s}  \underbrace{\left(\bm{n} \otimes \bm{s} \right) \cddot \mathbb{E}_{\kappa}}_{M_{ns}} \right\}  \, \d{s} \\
&= \sum_{\Gamma_{i} \subset \partial \Omega} \left[ M_{ns} v_w \right]_{\partial \Gamma_i} + \int_{\partial \Omega} \left\{ \diffp{v_w}{n} M_{nn}  - v_w \, \diffp{M_{ns}}{s} \right\} \, \d{s}
\end{aligned}
\end{equation}
Combining equations \eqref{eq:line1_wf1_ten}, \eqref{eq:line1_wf2_ten} and \eqref{eq:line1_ipbc_ten} the final expression which makes appear the dynamic boundary terms (forces and momenta) is found:
\begin{equation}
\int_{\Omega} v_w \diffp{\alpha_w}{t} \, \d\Omega  =  -  \int_{\Omega} \mathrm{Grad}(\mathrm{grad}(v_w))  \cddot \mathbb{E}_{\kappa} \, \d\Omega  +  \int_{\partial \Omega} \left\{ \diffp{v_w}{n} M_{nn}  + v_w \, \widetilde{q}_n \right\}  \, \d{s}  + \sum_{\Gamma_{i} \subset \partial \Omega} \left[ M_{ns} v_w \right]_{\partial \Gamma_i}.
\end{equation}
If the boundary is regular, the final expression simplifies:  
\begin{equation}
\int_{\Omega} v_w \diffp{\alpha_w}{t}  \, \d\Omega =  -  \int_{\Omega} \mathrm{Grad}(\mathrm{grad}(v_w))  \cddot \mathbb{E}_{\kappa} \, \d\Omega +  \int_{\partial \Omega} \left\{ \diffp{v_w}{n} M_{nn}  + v_w \, \widetilde{q}_n \right\} \, \d{s}. 
\end{equation}

So the final weak form obtained from system \eqref{eq:PH_sys_Kir_Ten} is written as:
\begin{equation}
\label{eq:WF_Kir_Dyn}
\begin{cases}
\displaystyle\int_{\Omega} v_w \diffp{\alpha_w}{t} \, \d\Omega  &=  -  \displaystyle\int_{\Omega} \mathrm{Grad}(\mathrm{grad}(v_w))  \cddot \mathbb{E}_{\kappa} \, \d\Omega +  \displaystyle\int_{\partial \Omega} \left\{ \diffp{v_w}{n} M_{nn}  + v_w \, \widetilde{q}_n \right\}   \, \d{s},  \vspace{2mm}\\
\displaystyle\int_{\Omega} \mathbb{V}_{\kappa} \cddot \diffp{\mathbb{A}_{\kappa}}{t} \, \d\Omega &= \displaystyle\int_{\Omega} \mathbb{V}_{\kappa} \cddot \mathrm{Grad}(\mathrm{grad}(e_w)) \, \d\Omega. 
\end{cases}
\end{equation}

The control inputs $\bm{u}_\partial$ and the corresponding conjugate outputs $\bm{y}_\partial$ are:
\[\bm{u}_\partial = 
\begin{pmatrix}
\widetilde{q}_n \\
M_{nn} \\
\end{pmatrix}_{\partial \Omega}, \qquad
\bm{y}_\partial = 
\begin{pmatrix}
w_t \\
\displaystyle \diffp{w_t}{n} \\
\end{pmatrix}_{\partial \Omega}.
\]

\subsubsection{Boundary control through kinematic variables}
\label{sec:BC_displ}
Alternatively, the same procedure can be performed on the second line of the system to make appear the kinematic boundary conditions, i.e. the value of the vertical velocity and its normal derivative along the border. Once the necessary calculations are carried out, the following result is found: 
\begin{equation}
\label{eq:WF_Kir_Kin}
\begin{cases}
\displaystyle\int_{\Omega} v_w \diffp{\alpha_w}{t} \, \d\Omega &=  \displaystyle\int_{\Omega} -v_w \, \mathrm{div}(\mathrm{Div}(\mathbb{E}_{\kappa})) \,  \d\Omega, \vspace{2mm}\\
\displaystyle\int_{\Omega} \mathbb{V}_{\kappa} \cddot \diffp{\mathbb{A}_{\kappa}}{t}   \d\Omega &= \displaystyle\int_{\Omega} \mathrm{div}(\mathrm{Div}(\mathbb{V}_{\kappa})) \; e_w  d\Omega +  \displaystyle\int_{\partial \Omega} \left\{ {v}_{M_{nn}} \diffp{w_t}{n}  + v_{\widetilde{q}_n} w_t \right\} \  \d{s}. 
\end{cases}
\end{equation}
where $v_{M_{nn}} = \left(\bm{n} \otimes \bm{n} \right) \cddot \mathbb{V}_{\kappa}  \; $ and $ \; v_{\widetilde{q}_n} = - \displaystyle \mathrm{Div}(\mathbb{V}_{\kappa}) \cdot \bm{n} - \diffp{ v_{M_{ns}} }{s}$ with $v_{M_{ns}} = \left(\bm{n} \otimes \bm{s} \right) \cddot \mathbb{E}_{\kappa}$. The control inputs $\bm{u}_\partial$ and the corresponding conjugate outputs $\bm{y}_\partial$ are: 
\[\bm{u}_\partial = 
\begin{pmatrix}
w_t \\
\displaystyle \diffp{w_t}{n} \\
\end{pmatrix}_{\partial \Omega}, \qquad
\bm{y}_\partial = 
\begin{pmatrix}
\widetilde{q}_n \\
M_{nn} \\
\end{pmatrix}_{\partial \Omega}.
\]

\subsection{Finite-dimensional port-Hamiltonian system}

{
	In this section, the discretization procedure is applied to formulation \eqref{eq:WF_Kir_Dyn}. The same procedure may be performed using formulation \eqref{eq:WF_Kir_Kin}. In Section \ref{sec:Eigen} both strategies will be used to compute the eigenvalues of a square plate. \\
	Test and co-energy variables are discretized using the same basis functions (Galerkin Method):
	\begin{equation}
	\begin{aligned}
	v_w &= \sum_{i = 1}^{N_w} \phi_w^i(x,y) \, v_w^i, \\
	\mathbb{V}_\kappa &= \sum_{i = 1}^{N_\kappa} \bm\Phi_\kappa^i(x,y) \, v_\kappa^i,\\
	\end{aligned} \qquad \quad
	\begin{aligned}
	e_w &= \sum_{i = 1}^{N_w} \phi_w^i(x,y) \, e_w^i(t), \\
	\mathbb{E}_\kappa &= \sum_{i = 1}^{N_\kappa} \bm\Phi_\kappa^i(x,y) \, e_\kappa^i(t),\\
	\end{aligned}
	\end{equation}
	The basis functions $\phi_w^i, \, \bm\Phi_\kappa^i, $ have to be chosen in a suitable function space $\mathcal{V}^h$ in the domain of operator $J$, i.e. $\mathcal{V}^h \subset \mathcal{V} \in \mathcal{D}(J)$. This will be discussed in Section \ref{sec:Num}. The discretized skew-symmetric bilinear form on the right side of \eqref{eq:WF_Kir_Dyn} then yields:
	\begin{equation}
	\bm{J}_d = 
	\begin{bmatrix}
	0 & -\bm{D}_{\mathrm{H}}^T \vspace{.3mm}\\
	\bm{D}_{\mathrm{H}} & 0 \vspace{.3mm}\\
	\end{bmatrix}.
	\end{equation}
	Matrix $\bm{D}_{\mathrm{H}}$ is computed in the following way:
	\begin{equation}
	\bm{D}_{\mathrm{H}}(i,j) = \int_{\Omega} \bm{\Phi}_{\kappa}^i : \mathrm{Grad}(\mathrm{grad}(\phi_w^j)) \d\Omega, \quad \in \mathbb{R}^{N_\kappa \times N_w},
	\end{equation}
	where the notation $A(i,j)$ indicates the entry in the matrix corresponding to the $i \, {\text{th}}$ row and $j \,{\text{th}}$ column. The energy variables are deduced from the co-energy variables: 
	\begin{equation}
	\alpha_w = \mu e_w, \qquad \mathbb{A}_{\kappa} = \mathbb{D}^{-1} \mathbb{E}_{\kappa},
	\end{equation}
	where $\mathbb{D}_{ijkl}$ is the symmetric bending rigidity tensor, the tensorial analogous of  matrix $\bm{D}$ defined in \eqref{eq:BenRigVec}. The symmetric bilinear form on the left side of \eqref{eq:WF_Kir_Dyn} becomes: 
	\begin{equation}
	\begin{gathered}
	\bm{M} = \text{diag}[\bm{M}_w,\, \bm{M}_\kappa], \quad \text{with} \\
	\begin{aligned}
	&\bm{M}_w(i,j) = \int_{\Omega} \mu \, \bm{\phi}_w^i \, \bm{\phi}_w^j \d\Omega, \; \in \mathbb{R}^{N_w \times N_w}, \\
	&\bm{M}_\kappa(i,j) = \int_{\Omega}  \left( \mathbb{D}^{-1} \bm{\Phi}_\kappa^i \right) \cddot \bm{\Phi}_\kappa^j \d\Omega, \; \in \mathbb{R}^{N_\kappa \times N_\kappa}.\\
	\end{aligned}
	\end{gathered}
	\end{equation}
	The boundary variables are then discretized as:
	\begin{equation}
	\begin{aligned}
	\widetilde{q}_n = \sum_{i = 1}^{N_{\widetilde{q}_n}} \phi_{\widetilde{q}_n}^i(s) \; \widetilde{q}_n^i, \qquad
	M_{nn} = \sum_{i = 1}^{N_{M_{nn}}} \phi_{M_{nn}}^i(s) \; M_{nn}^i.
	\end{aligned}
	\end{equation}
	The variables are defined only over the boundary $\partial\Omega$. Consequently, the input matrix reads:
	\begin{equation}
	\bm{B} = \begin{bmatrix}
	\bm{B}_{\widetilde{q}_n} & \bm{B}_{M_{nn}}\\
	0 & 0 \\
	\end{bmatrix}. 
	\end{equation}
	The inner components are computed as:
	\begin{equation}
	\begin{aligned}
	\bm{B}_{\widetilde{q}_n}(i,j) &= \int_{\partial\Omega} {\phi}_{w}^i \, {\phi}_{q_n}^j \d{s}, \quad \in \mathbb{R}^{N_w \times N_{q_n}} , \\
	\bm{B}_{M_{nn}}(i,j) &= \int_{\partial\Omega} \diffp{\phi_w^i}{n} \, \phi_{M_{nn}}^j \d{s}, \quad \in \mathbb{R}^{N_w \times N_{M_{nn}}} . \\
	\end{aligned}
	\end{equation}
	The final port-Hamiltonian system, as defined in \cite{beattie2018linear} is written as:
	\begin{equation}
	\label{eq:PHdiscr_Kir}
	\begin{aligned}
	\bm{M} \dot{\bm{e}} &= \bm{J}_d  \,\bm{e} + \bm{B} \, \bm{u}_{\partial}, \\
	\bm{y}_{\partial} &= \bm{B}^T \, \bm{e},
	\end{aligned} 
	\end{equation}
	where $\bm{e} = \left(e_w^1, \cdots, e_{\kappa}^{N_{\kappa}}\right)^T$ and $\bm{u}_{\partial} = \left(\widetilde{q}_n^1, \dots, M_{nn}^{N_{M_{nn}}}\right)^T$ are the concatenations of the degrees of freedom for the different variables. The discrete Hamiltonian is then found as:
	\begin{equation}
	\label{eq:Hd_Kir}
	\begin{aligned}
	H_d &= \frac{1}{2} \int_{\Omega} \left\{ \alpha_w e_w + \mathbb{A}_{\kappa} \cddot \mathbb{E}_\kappa \right\} \d\Omega \\
	&=  \frac{1}{2}  \left\{ \bm{e}_w^T \, \bm{M}_w \, \bm{e}_w + \bm{e}_\kappa^T \, \bm{M}_\kappa \, \bm{e}_\kappa \right\} \\
	&=  \frac{1}{2} \, \bm{e}^T  \bm{M}  \bm{e}.
	\end{aligned}
	\end{equation}
}
Using equations \eqref{eq:PHdiscr_Kir}, \eqref{eq:Hd_Kir} the time derivative of the Hamiltonian is given by the scalar product of the boundary flows:
\begin{equation}
\dot{H}_d = \bm{y}_{\partial}^T \bm{u}_{\partial}.
\end{equation}
The above Equation is equivalent to the energy balance of the continuous system, expressed by \eqref{eq:energyBal_Kir}. Definition \eqref{eq:Hd_Kir}, together with system \eqref{eq:PHdiscr_Kir} are the finite-dimensional equivalent of \eqref{eq:H_Kir_Ten} and  \eqref{eq:PH_sys_Kir_Ten}. The discretized system obtained via PFEM shares the {port-Hamiltonian structure} of the original infinite-dimensional system, the discretization method is therefore structure preserving. 

\begin{figure}[p]%
	\minipage{0.25\textwidth}%
	\includegraphics[width=\linewidth]{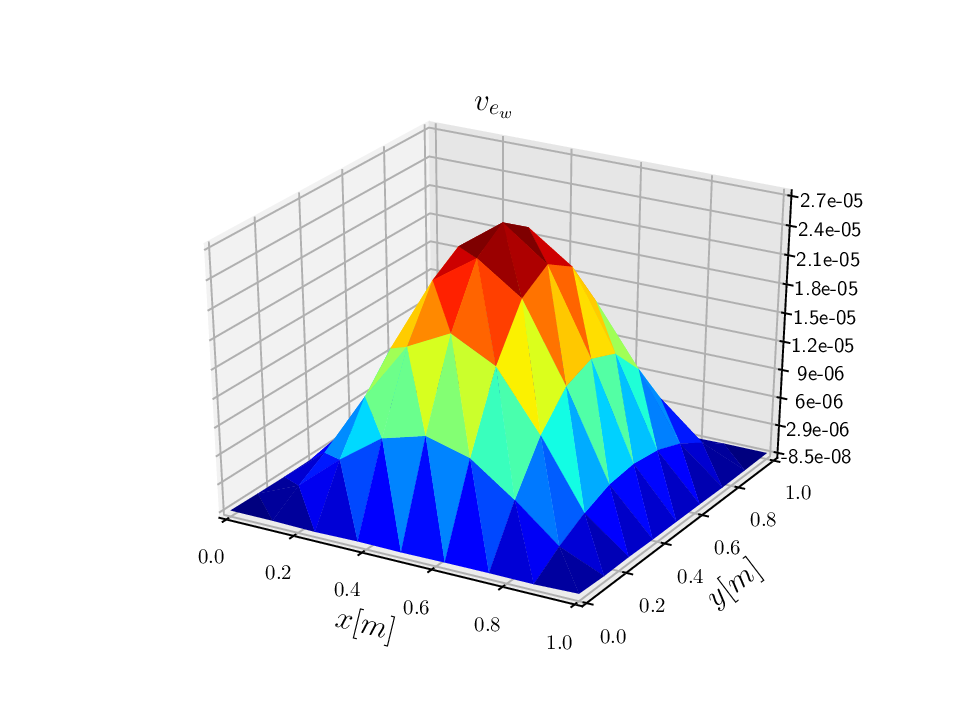}
	\caption*{$\widehat{\omega}_{1}$}\label{fig:CSCS1}%
	\endminipage
	\minipage{0.25\textwidth}%
	\includegraphics[width=\linewidth]{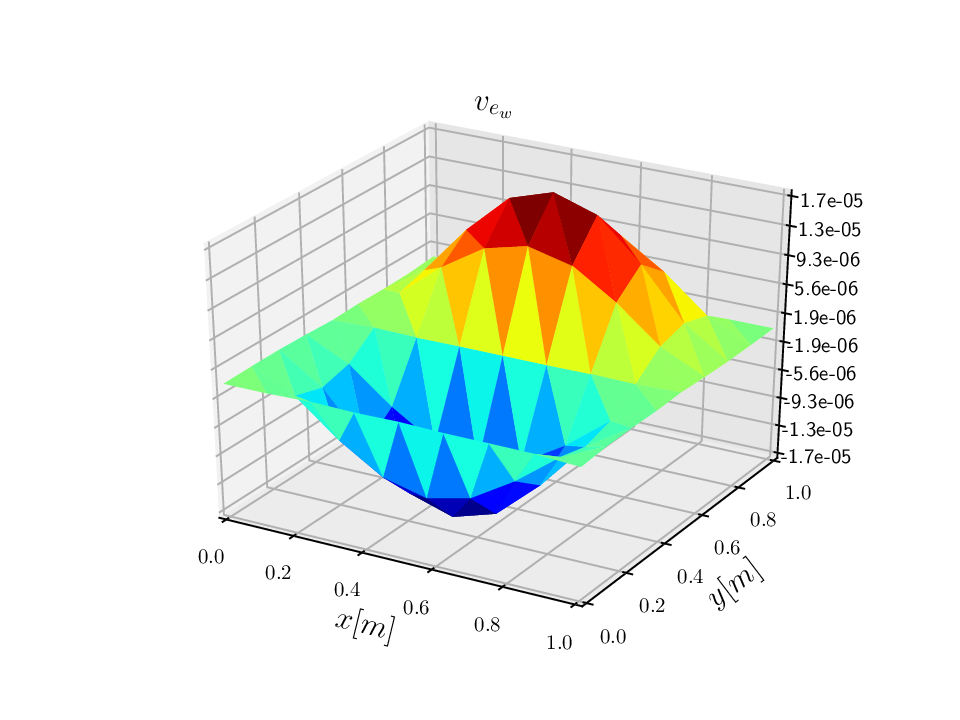}
	\caption*{$\widehat{\omega}_{2}$}\label{fig:CSCS2}%
	\endminipage
	\minipage{0.25\textwidth}%
	\includegraphics[width=\linewidth]{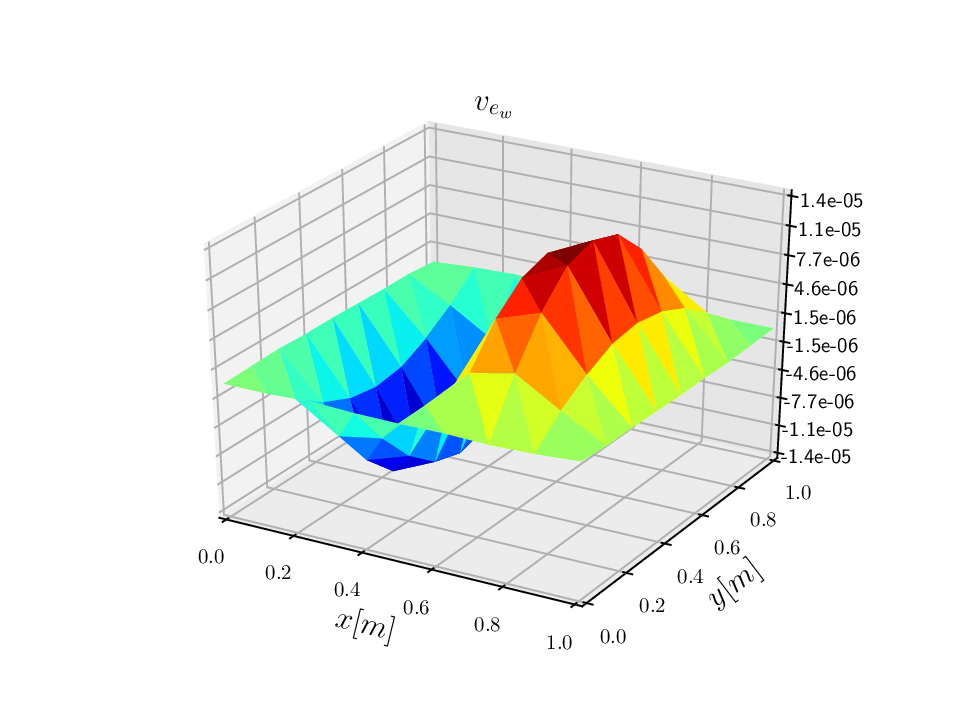}
	\caption*{$\widehat{\omega}_{3}$}\label{fig:CSCS3}%
	\endminipage 
	\minipage{0.25\textwidth}%
	\includegraphics[width=\linewidth]{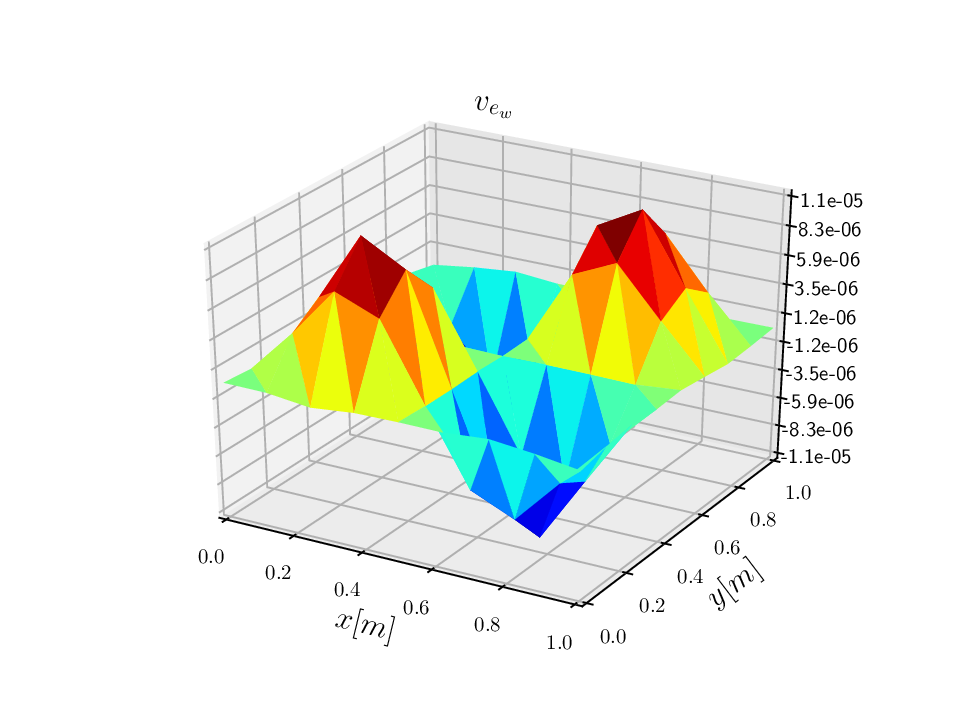}
	\caption*{$\widehat{\omega}_{4}$}\label{fig:CSCS4}%
	\endminipage
	\caption[Eigenvectors for CSCS]{Eigenvectors for the CSCS case.}%
	\label{fig:CSCS}%
\end{figure}

\begin{figure}[p]%
	\minipage{0.25\textwidth}%
	\includegraphics[width=\linewidth]{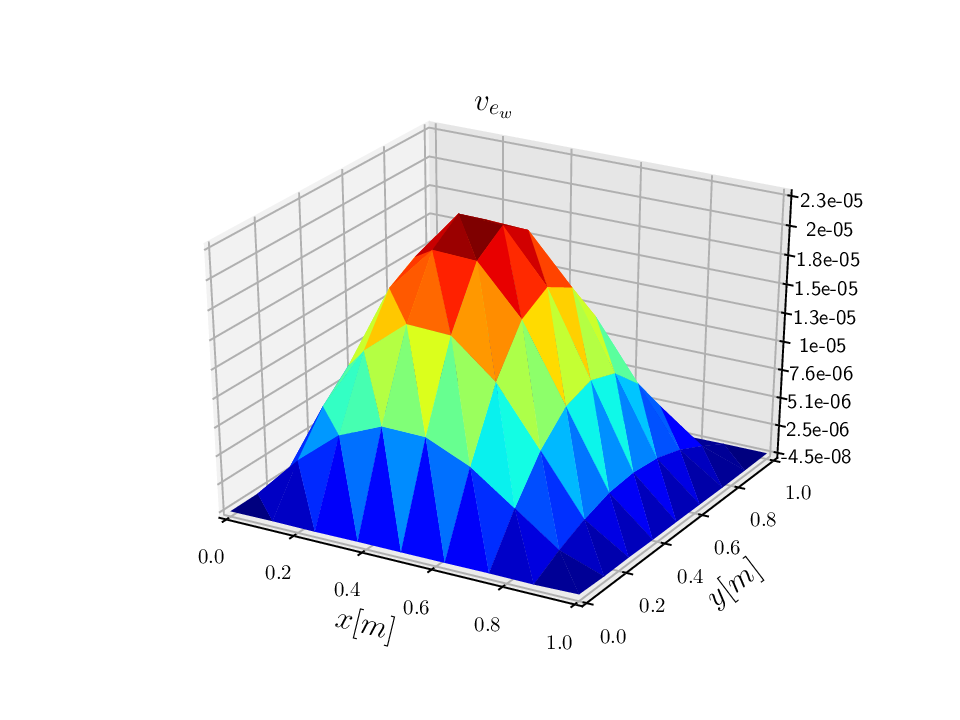}
	\caption*{$\widehat{\omega}_{1}$}\label{fig:SSCS1}%
	\endminipage
	\minipage{0.25\textwidth}%
	\includegraphics[width=\linewidth]{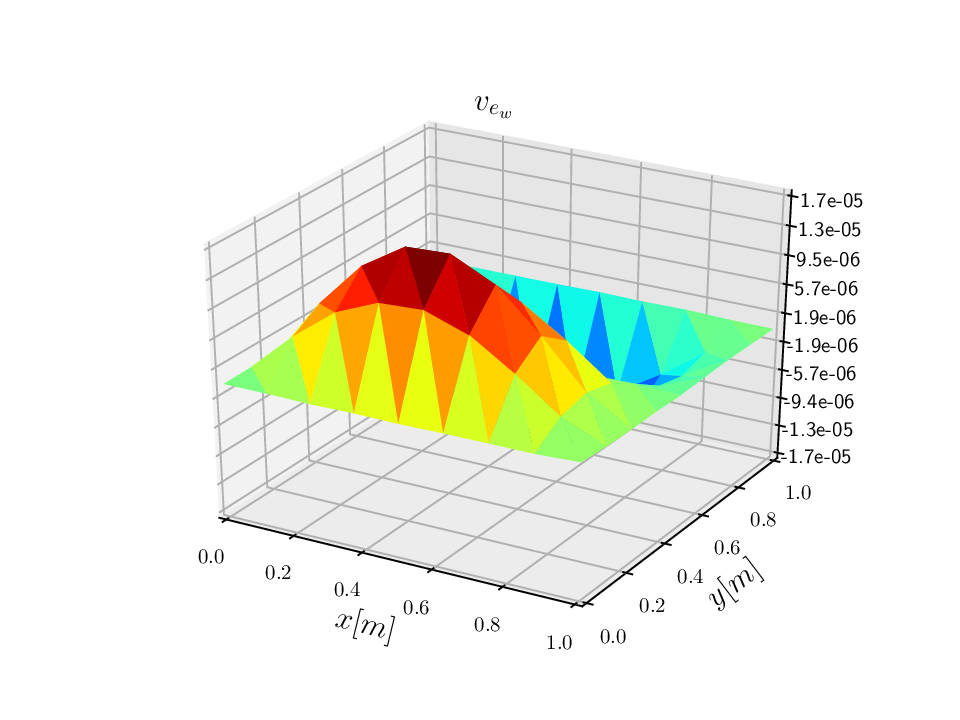}
	\caption*{$\widehat{\omega}_{2}$}\label{fig:SSCS2}%
	\endminipage
	\minipage{0.25\textwidth}%
	\includegraphics[width=\linewidth]{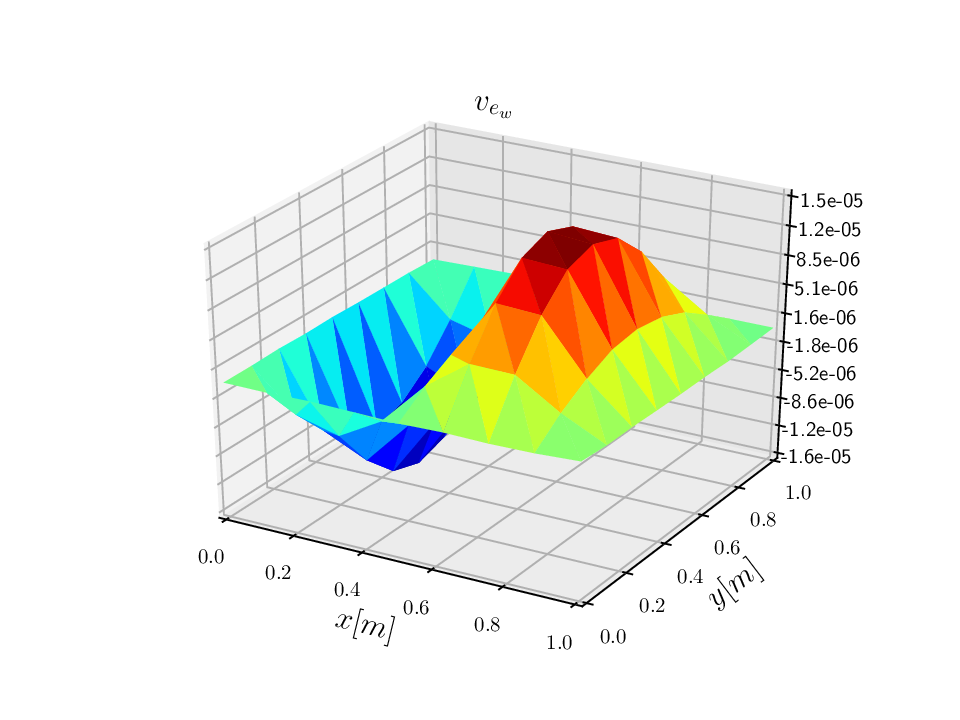}
	\caption*{$\widehat{\omega}_{3}$}\label{fig:SSCS3}%
	\endminipage 
	\minipage{0.25\textwidth}%
	\includegraphics[width=\linewidth]{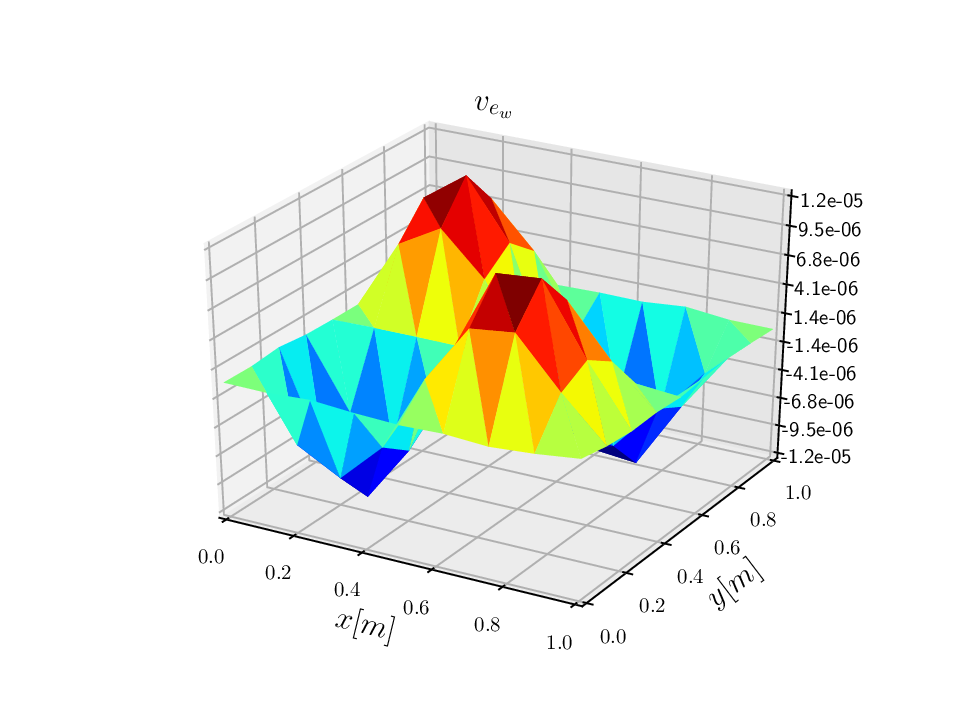}
	\caption*{$\widehat{\omega}_{4}$}\label{fig:SSCS4}%
	\endminipage
	\caption[Eigenvectors for SSCS]{Eigenvectors for the SSCS case.}%
	\label{fig:SSCS}%
\end{figure}

\begin{figure}[p]%
	\minipage{0.25\textwidth}%
	\includegraphics[width=\linewidth]{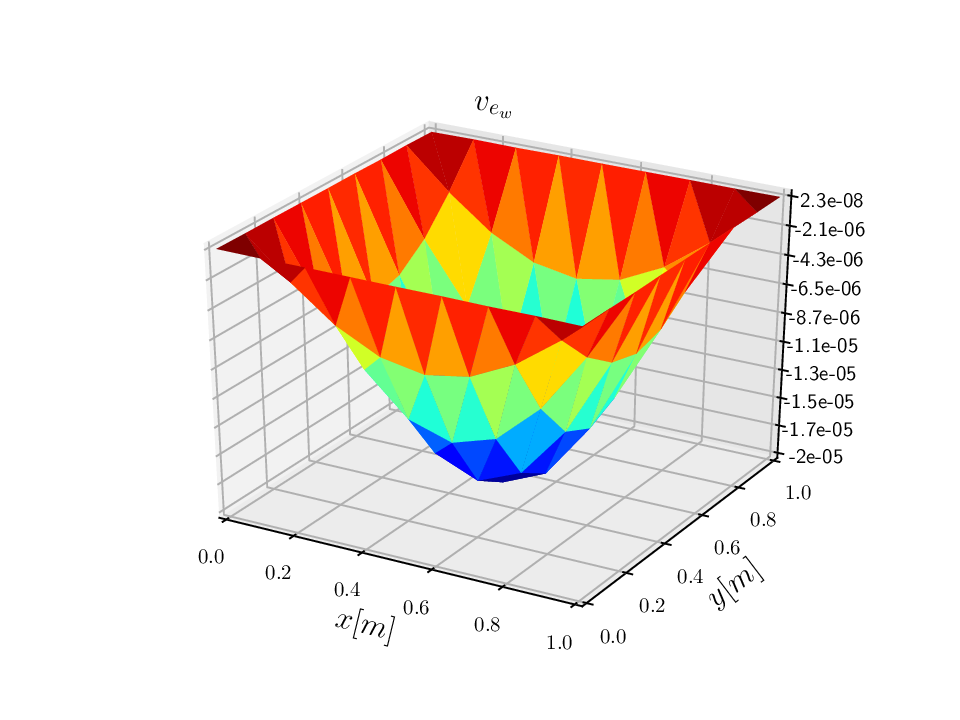}
	\caption*{$\widehat{\omega}_{1}$}\label{fig:SSSS1}%
	\endminipage
	\minipage{0.25\textwidth}%
	\includegraphics[width=\linewidth]{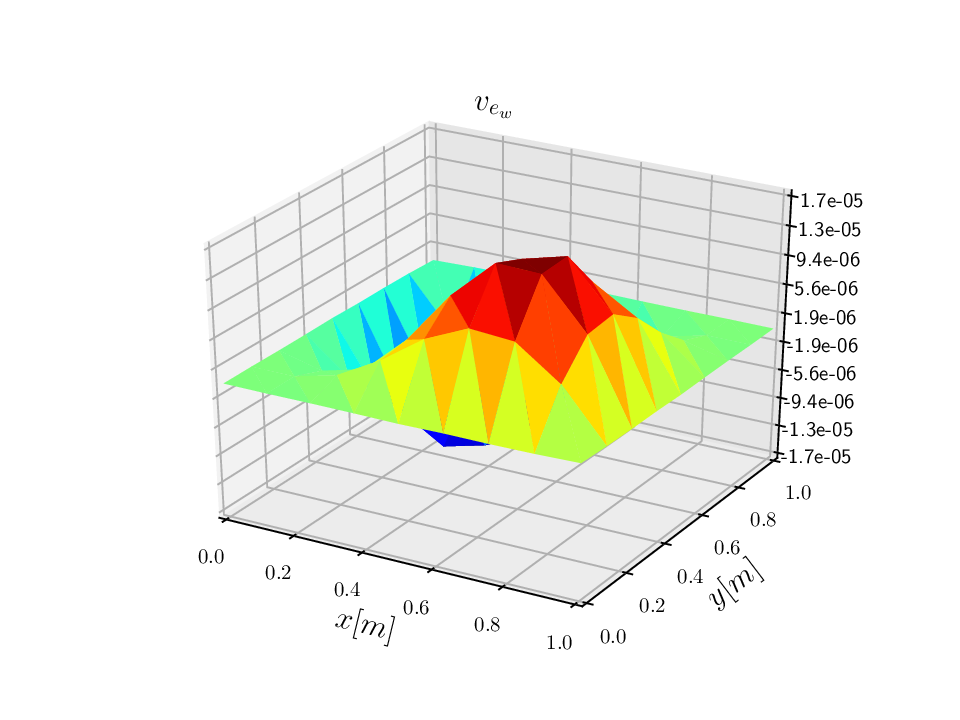}
	\caption*{$\widehat{\omega}_{2}$}\label{fig:SSSS2}%
	\endminipage
	\minipage{0.25\textwidth}%
	\includegraphics[width=\linewidth]{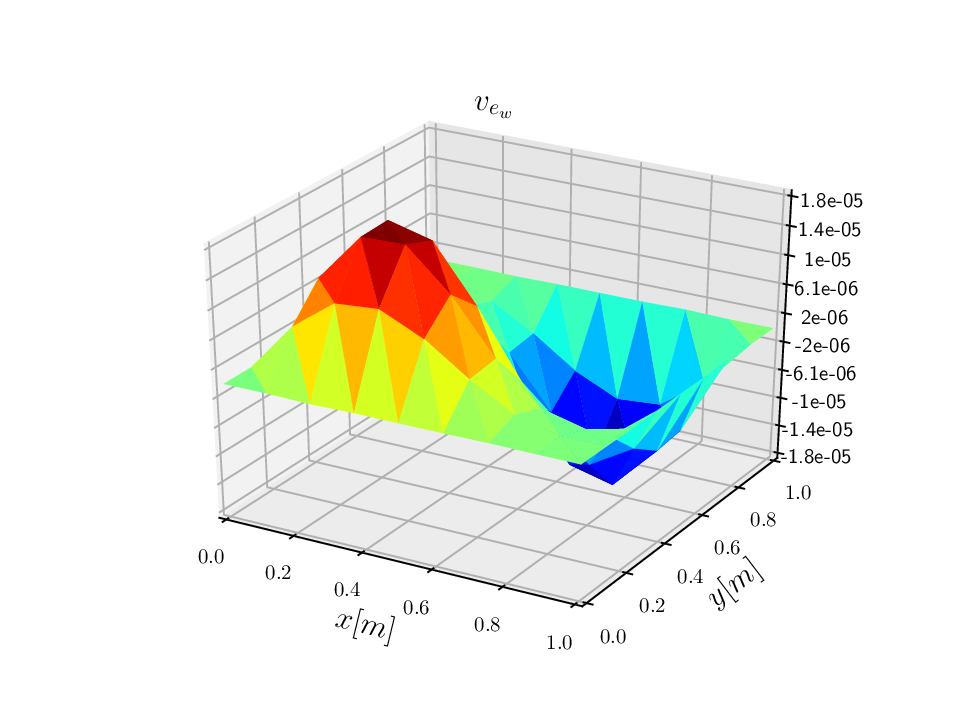}
	\caption*{$\widehat{\omega}_{3}$}\label{fig:SSSS3}%
	\endminipage 
	\minipage{0.25\textwidth}%
	\includegraphics[width=\linewidth]{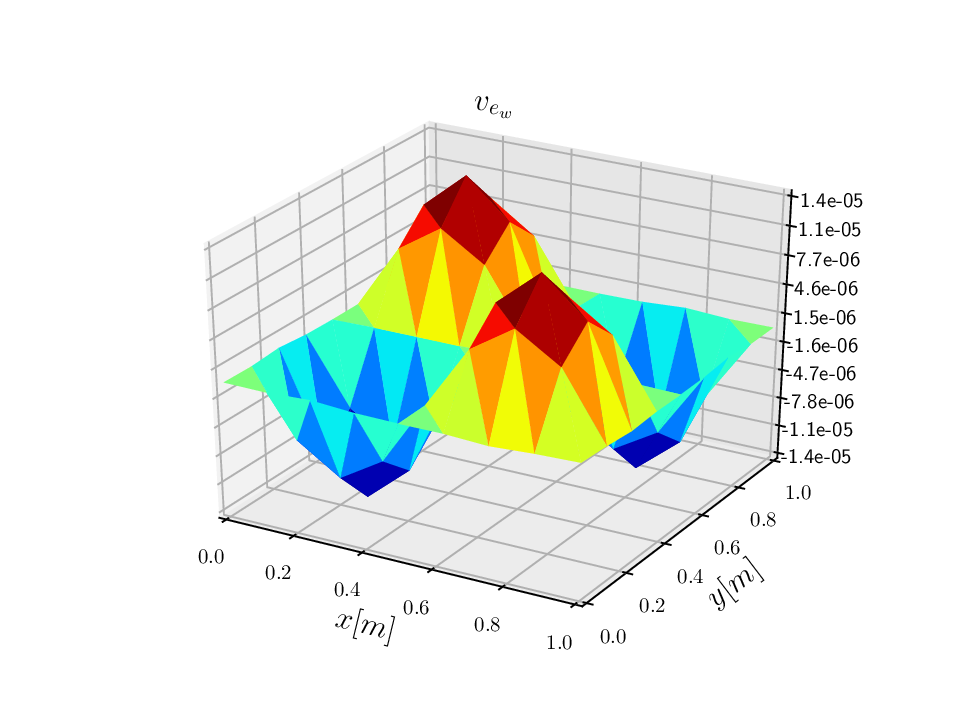}
	\caption*{$\widehat{\omega}_{4}$}\label{fig:SSSS4}%
	\endminipage
	\caption[Eigenvectors for SSSS]{Eigenvectors for the SSSS case.}%
	\label{fig:SSSS}%
\end{figure}

\begin{figure}[p]%
	\minipage{0.25\textwidth}%
	\includegraphics[width=\linewidth]{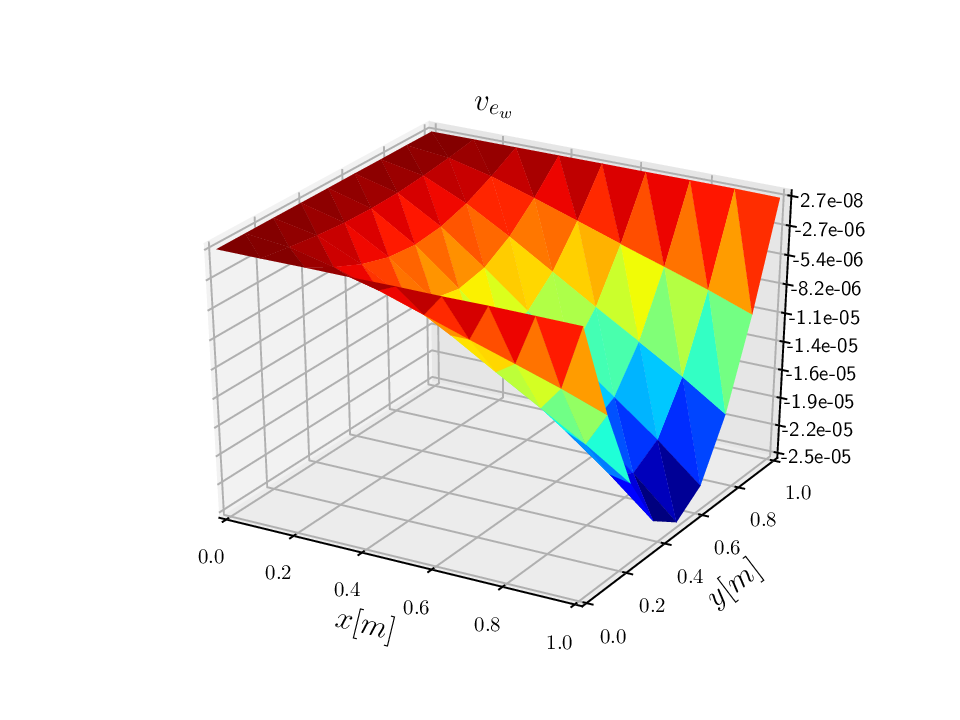}
	\caption*{$\widehat{\omega}_{1}$}\label{fig:CSFS1}%
	\endminipage
	\minipage{0.25\textwidth}%
	\includegraphics[width=\linewidth]{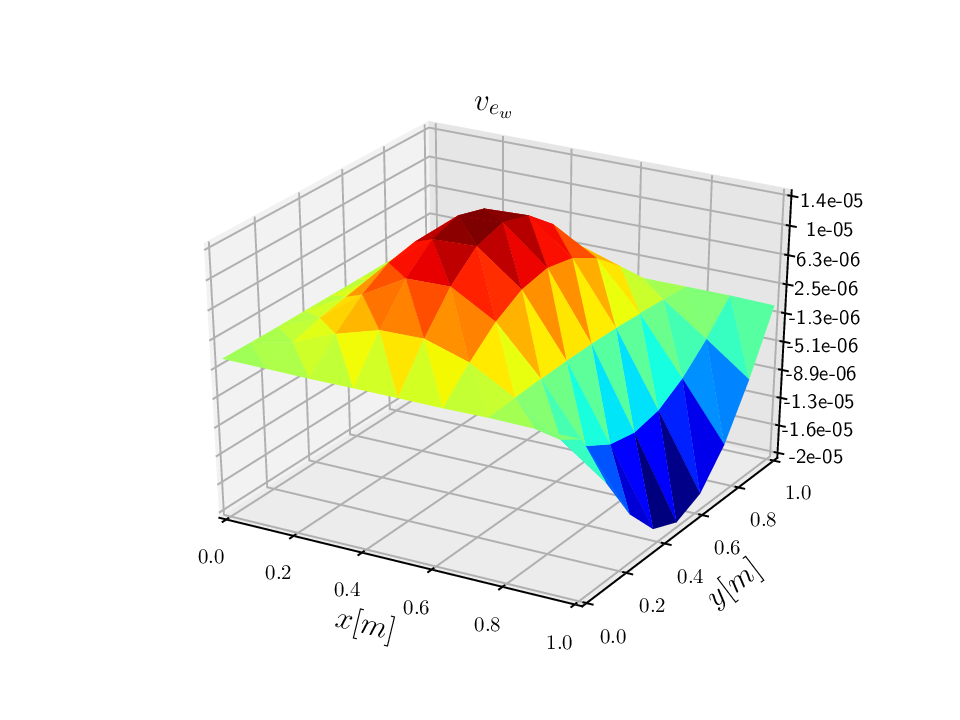}
	\caption*{$\widehat{\omega}_{2}$}\label{fig:CSFS2}%
	\endminipage
	\minipage{0.25\textwidth}%
	\includegraphics[width=\linewidth]{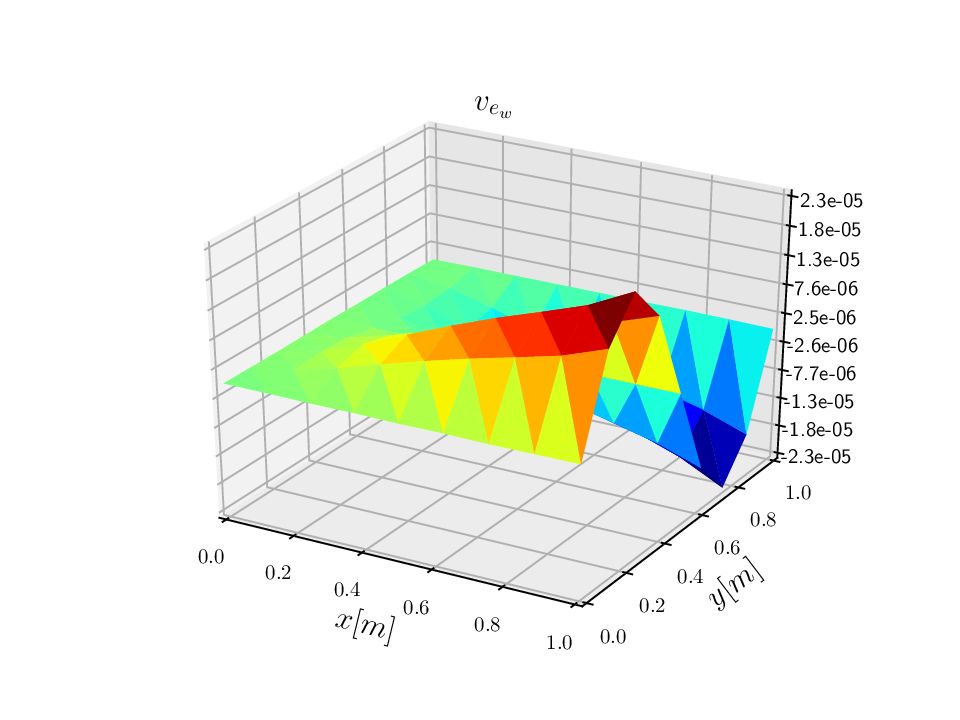}
	\caption*{$\widehat{\omega}_{3}$}\label{fig:CSFS3}%
	\endminipage 
	\minipage{0.25\textwidth}%
	\includegraphics[width=\linewidth]{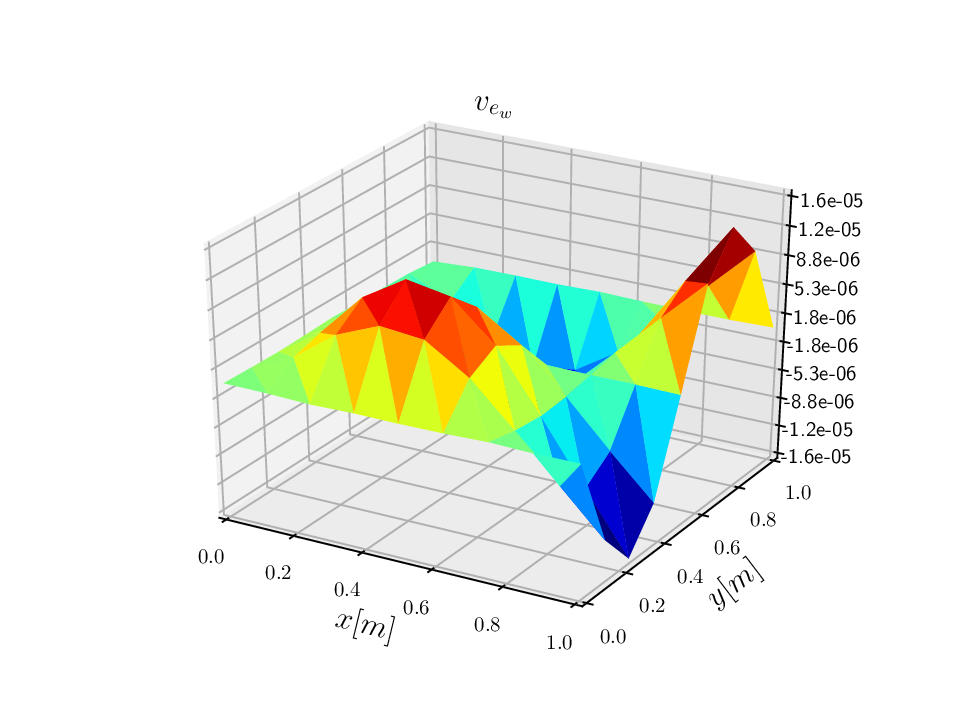}
	\caption*{$\widehat{\omega}_{4}$}\label{fig:CSFS4}%
	\endminipage
	\caption[Eigenvectors for CSFS]{Eigenvectors for the CSFS case.}%
	\label{fig:CSFS}%
\end{figure}

\begin{figure}[ht]%
	\minipage{0.25\textwidth}%
	\includegraphics[width=\linewidth]{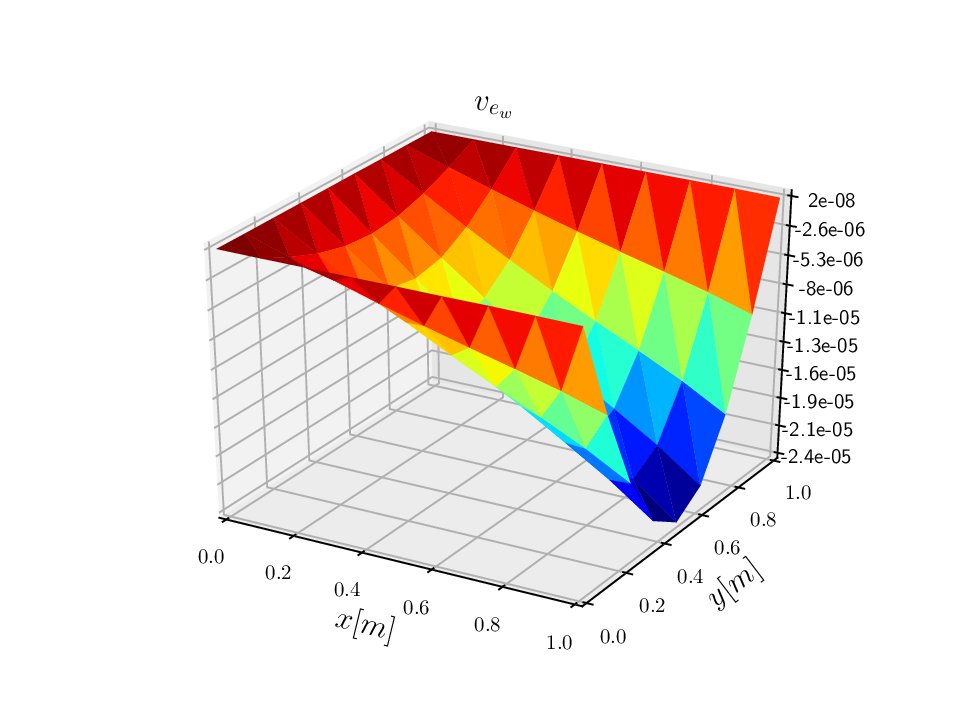}
	\caption*{$\widehat{\omega}_{1}$}\label{fig:SSFS1}%
	\endminipage
	\minipage{0.25\textwidth}%
	\includegraphics[width=\linewidth]{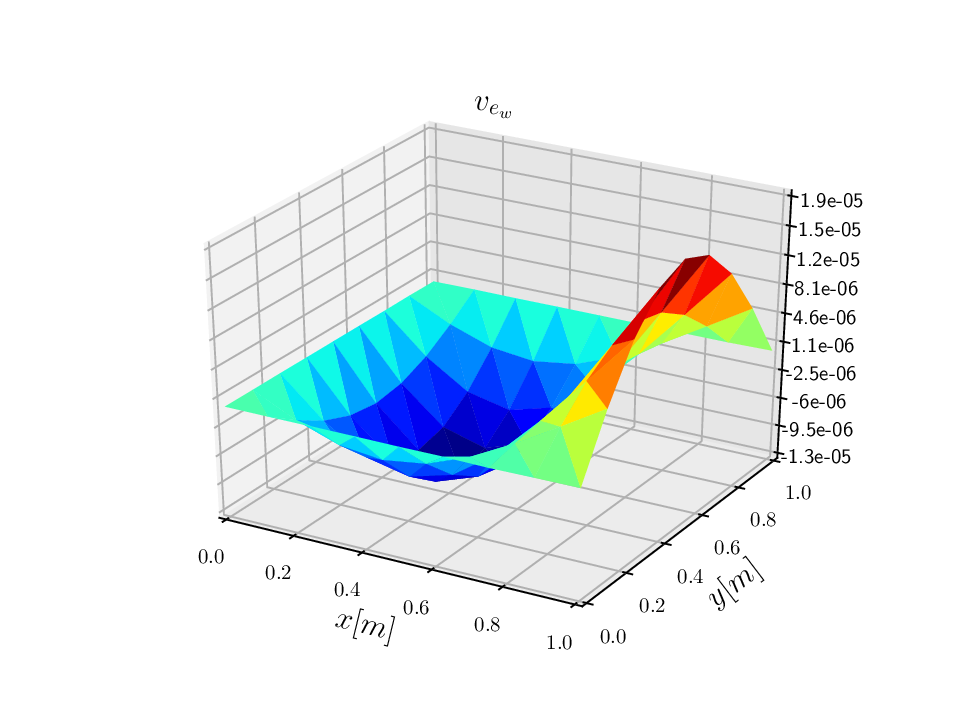}
	\caption*{$\widehat{\omega}_{2}$}\label{fig:SSFS2}%
	\endminipage
	\minipage{0.25\textwidth}%
	\includegraphics[width=\linewidth]{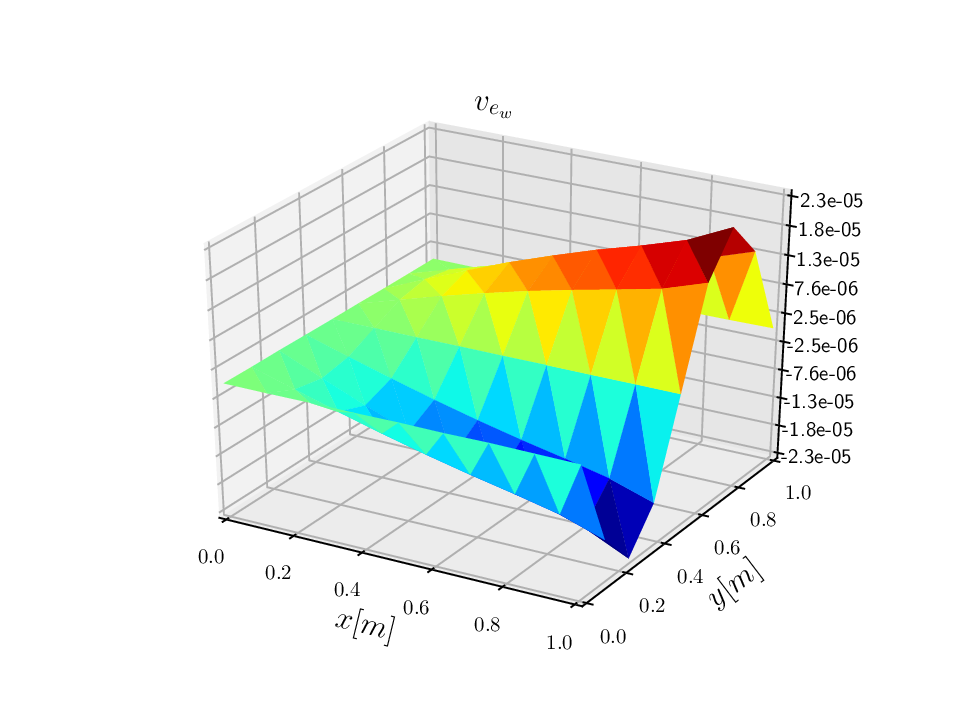}
	\caption*{$\widehat{\omega}_{3}$}\label{fig:SSFS3}%
	\endminipage 
	\minipage{0.25\textwidth}%
	\includegraphics[width=\linewidth]{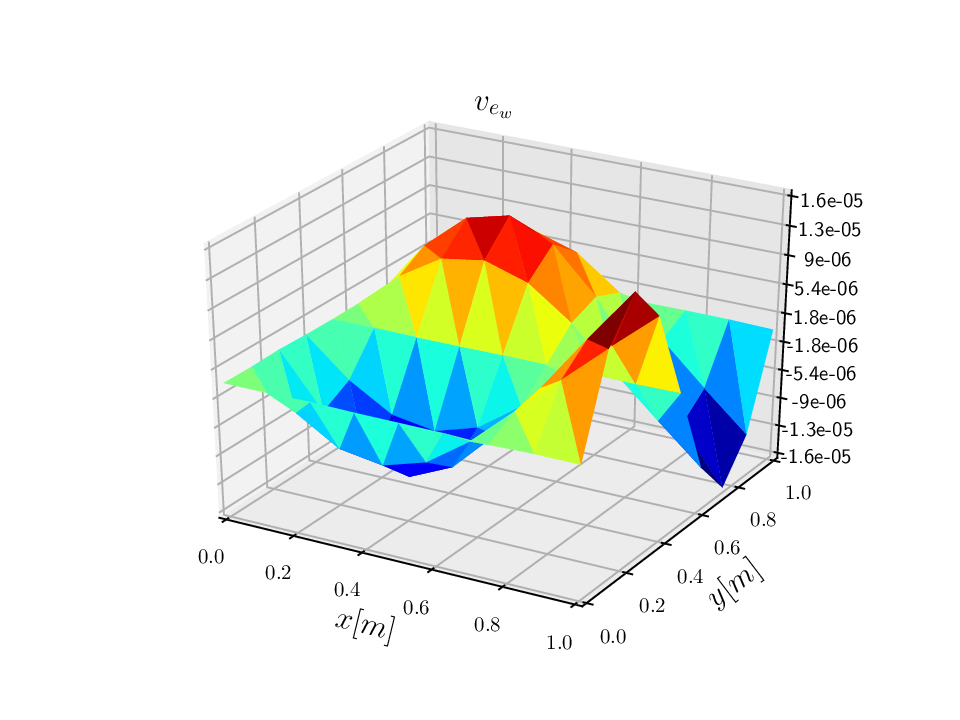}
	\caption*{$\widehat{\omega}_{4}$}\label{fig:SSFS4}%
	\endminipage
	\caption[Eigenvectors for SSFS]{Eigenvectors for the SSFS case.}%
	\label{fig:SSFS}%
\end{figure}

\begin{figure}[h]%
	\minipage{0.25\textwidth}%
	\includegraphics[width=\linewidth]{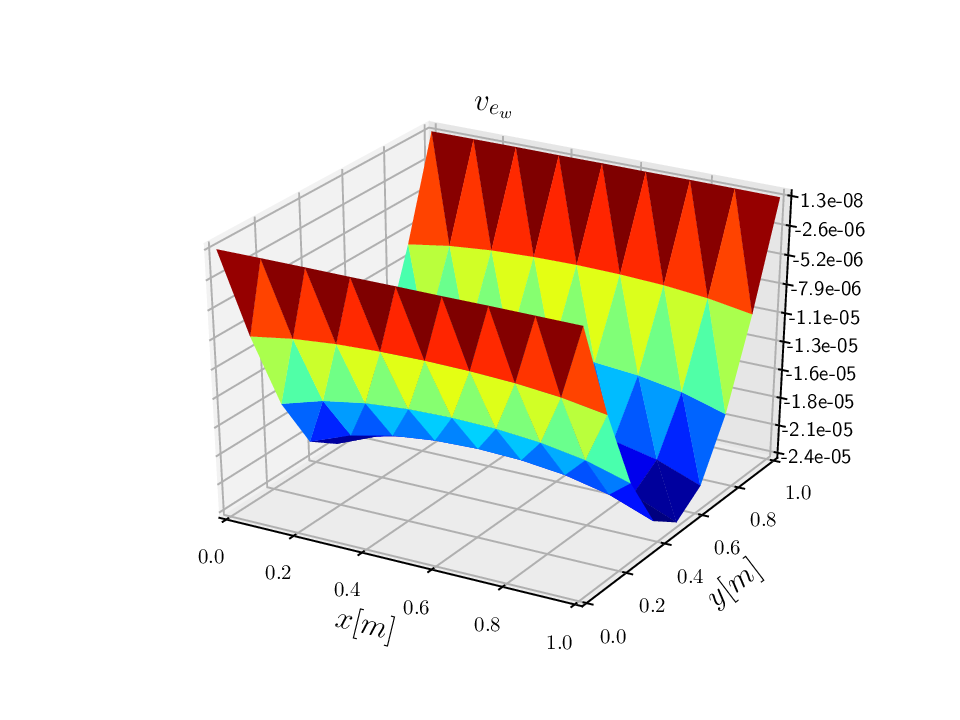}
	\caption*{$\widehat{\omega}_{1}$}\label{fig:FSFS1}%
	\endminipage
	\minipage{0.25\textwidth}%
	\includegraphics[width=\linewidth]{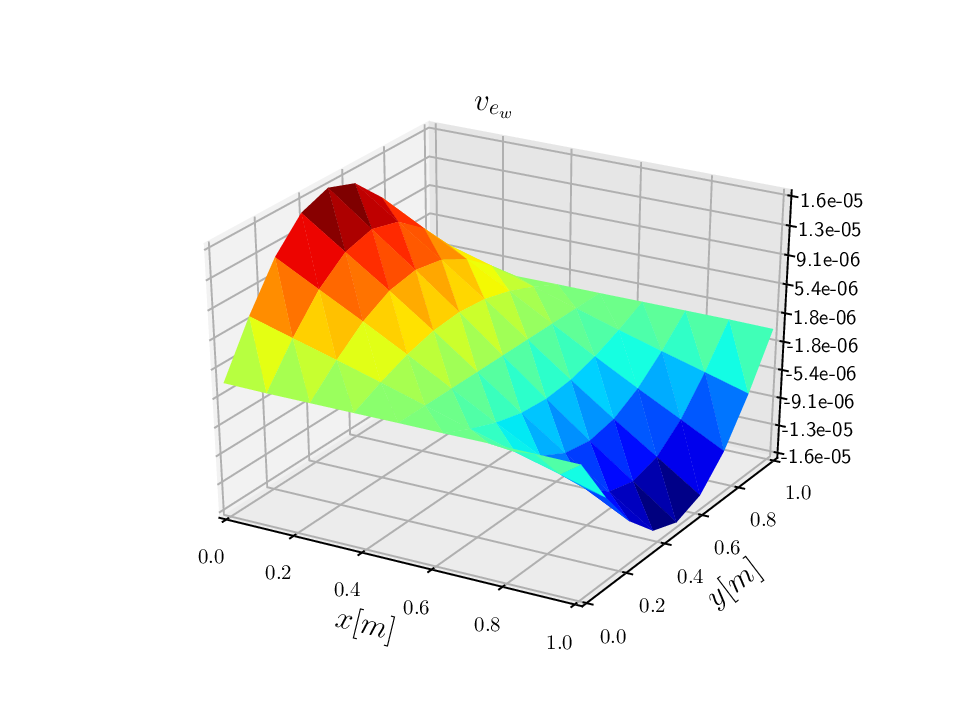}
	\caption*{$\widehat{\omega}_{2}$}\label{fig:FSFS2}%
	\endminipage
	\minipage{0.25\textwidth}%
	\includegraphics[width=\linewidth]{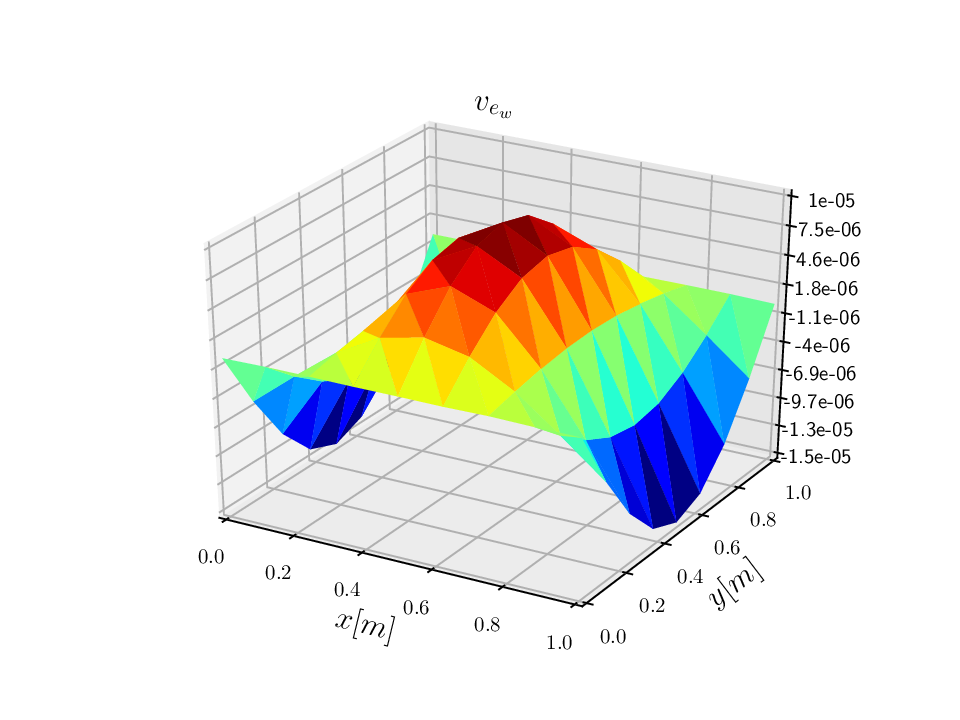}
	\caption*{$\widehat{\omega}_{3}$}\label{fig:FSFS3}%
	\endminipage 
	\minipage{0.25\textwidth}%
	\includegraphics[width=\linewidth]{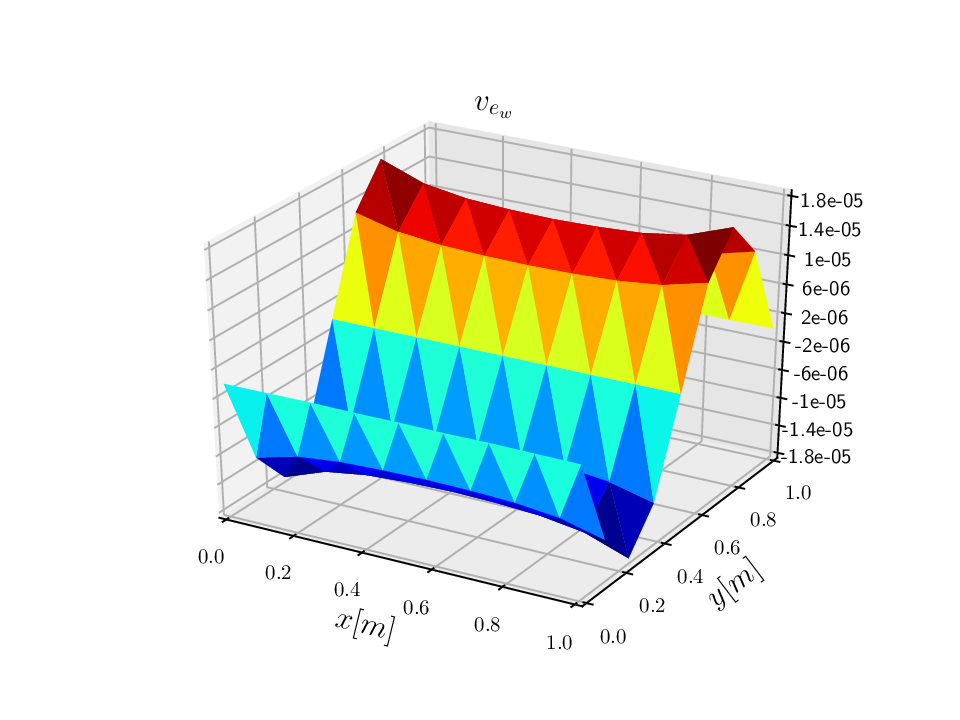}
	\caption*{$\widehat{\omega}_{4}$}\label{fig:FSFS4}%
	\endminipage
	\caption[Eigenvectors for FSFS]{Eigenvectors for the FSFS case.}%
	\label{fig:FSFS}%
\end{figure}

{
	\section{Numerical studies}
	\label{sec:Num}
	In this section we illustrate numerically the consistency of discrete model obtained with PFEM. For this purpose computation of the eigenvalues of a square plate and time-domain simulations for several boundary conditions are presented. 
	\subsection{Finite Element Choice}
	The domain of the operator $J$ in \eqref{eq:PH_sys_Kir_Ten} is $\mathcal{D}(J) = H^{2}(\Omega) \times  H^{\text{div Div}}(\Omega, \mathbb{R}^{2 \times 2}_{\text{sym}})$ and boundary conditions. 
	\begin{remark}
		It has to be appointed that, to the best of authors' knowledge, the space $H^{\text{div Div}}(\Omega, \mathbb{R}^{2 \times 2}_{\text{sym}})$ has never addressed in the mathematical literature. For this reason $H^{2}(\Omega)$ conforming finite elements were used to deal with this problem numerically. 
	\end{remark}
	A suitable choice for the functional space is thus:
	\begin{equation}
	(v_w, \,\mathbb{V}_\kappa) \in H^{2}(\Omega) \times H^{2}(\Omega, \mathbb{R}^{2 \times 2}_{\text{sym}}) \equiv \mathscr{H},
	\end{equation}
	since $\mathscr{H} \subset \mathcal{D}(J)$. 
	The $H^2$ conforming finite elements (like the Hermite, Bell or Argyris finite elements) do not satisfy the proper equivalence properties to give a simple relationship between the reference basis and nodal basis on a general cell \cite{KirbyFE}. The Firedrake library \cite{firedrake} was used to implement the numerical analysis as it provides functionalities to automate the generalized mappings for these elements. \\ \\
	Then for the Finite Element choice, denote
	\[ H_r^k(\mathbb{P}_l, \Omega) = \{ v \in H^k(\Omega)|\; v_{|T} \in \mathbb{P}_l \; \forall T \in \mathbb{T}_r \} 
	\]
	the finite element space which is a subspace of $H^k(\Omega)$, based on the shape function space of piecewise polynomials of degree $l$. The shape function space is defined over the mesh $\mathbb{T}_r = \bigcup_i T_i$, where the cells $T_i$ are triangles. These spaces can be scalar-valued or symmetric matrix valued, depending on the variables to be discretized. The parameter $r$ is the average size of a mesh element.
	All the variables, i.e. the velocity $e_w$ and the momenta tensor $\mathbb{V}_\kappa$ as well as the corresponding test functions, are discretized by the same finite element space, the Bell finite element space \cite{Bell}, denoted $H_r^2(\mathbb{P}_5, \Omega)$. For this element the field is computed using quintic polynomials whose degrees of freedom are the values of the function, its gradient and its Hessian at the vertex of each triangular element. To deal with mixed boundary conditions Lagrange multipliers have to be introduced (the reader can refer to \cite{BrugnoliMin}, section 4.3 for an explanation). The multipliers are therefore discretized by using second degree Lagrange polynomials defined over the boundary $H_r^1(\mathbb{P}_2, \partial\Omega)$. 
	\subsection{Eigenvalues Computation}
	\label{sec:Eigen}
	The test case for this analysis is a simple square plate of side $L$, a benchmark problem which has been studied in \cite{LeissaRect, LimEigen} for different boundary conditions on each plate side. The possible cases are the following: 
	\begin{itemize}
		\item clamped side (C), for which $w_t = 0, \, \diffp{w_t}{n}=0$; 
		\item simply supported side (S), $w_t = 0, \, M_{nn}=0$;
		\item free side (F), $\widetilde{q}_n = 0, \, M_{nn}=0$.
	\end{itemize}
	In order to compare our results the eigenfrequencies $\omega_{n}^h$ are computed in the following non-dimensional form:
	\begin{equation}
	\widehat{\omega}_{n}^h = 4 L^2 \omega_{n}^h \left(\frac{\rho h}{D} \right)^{1/2},
	\end{equation}
	The only parameter which influences the results is the Poisson's ratio $\nu=0.3$. The reported non-dimensional frequencies are independent of the remaining geometrical and physical parameters.  The error is computed as:
	\begin{equation}
	\varepsilon = \frac{\text{abs}(\widehat{\omega}_{n}^h - \omega_{n}^{L})}{\omega_{n}^{L}},
	\end{equation}
	where $\omega_{n}^{L}$ are the eigenvalues computed in \cite{LeissaRect}. The results are computed either by using the forces and momenta as control \eqref{eq:WF_Kir_Dyn} or the vertical linear and angular velocity \eqref{eq:WF_Kir_Kin} (column Hessian and divDiv in Tables \ref{tab:eigNoFree}, \ref{tab:eigFree}).  The results are obtained using a regular mesh composed by 5 Bell element on each side. Hence, the state vector has a total dimension of 864. The dimension of the Lagrange multiplier vector depends on the boundary conditions upon consideration.  When using $H_r^2(\mathbb{P}_2, \partial\Omega)$ on the considered mesh, this number can vary from 0 to 80. The results obtained by using \eqref{eq:WF_Kir_Dyn} are in perfect agreement with the reference. This formulation was also used to compute the eigenvectors corresponding to the vertical velocity for the different cases under examination (see Figs \ref{fig:CSCS} to \ref{fig:FSFS}). For what concerns the weak formulation \eqref{eq:WF_Kir_Kin} the results deteriorate when a free condition (see Table \ref{tab:eigFree}) is present. 
}

\begin{table}[t]
	\centering	
	\begin{tabular}{|c|b g g|b g g|b g g|}
		\hline 
		
		\multirow{2}{*}{$n$} &   \multicolumn{3}{c|}{CSCS}	& \multicolumn{3}{c|}{SSCS} &	\multicolumn{3}{c|}{SSSS}	\\	
		&\cellcolor{white}Leissa & \cellcolor{white}Hessian & \cellcolor{white}divDiv  &\cellcolor{white}Leissa & \cellcolor{white}Hessian & \cellcolor{white}divDiv  &\cellcolor{white}Leissa & \cellcolor{white}Hessian & \cellcolor{white}divDiv \\
		\hline
		$\widehat{\omega}_{1}$&	28.946&	   28.950&	28.951&	23.646&	23.640&	23.646&	19.739&	19.730&	19.739\\
		$\widehat{\omega}_{2}$&	54.743&	   54.747&	54.744&	51.674&	51.666&	51.675&	49.348&	49.333&	49.348\\
		$\widehat{\omega}_{3}$&	69.32 &	   69.331&	69.330&	58.641&	58.641&	58.647&	49.348&	49.336&	49.348\\
		$\widehat{\omega}_{4}$&	94.584&	   94.602&	94.593&	86.126&	86.121&	86.138&	78.957&	78.920&	78.958\\
		$\widehat{\omega}_{5}$&	102.213&   102.245&	102.221&100.259&100.284&100.275&98.696&	98.692&	98.700\\
		$\widehat{\omega}_{6}$&	129.086&   129.141&	129.110&113.217&113.250&113.237&98.696&	98.709& 98.703\\		
		\hline 
	\end{tabular} 	
	\caption{Eigenvalues obtained with 5 Bell element per side for $\nu=0.3$, considering either the $\mathrm{Grad} \circ \mathrm{grad}$ formulation \eqref{eq:WF_Kir_Dyn}, either the $\mathrm{Div} \circ \mathrm{div}$ formulation \eqref{eq:WF_Kir_Kin}. For comparison reference \cite{LeissaRect} is considered.  \\
		\sqbox{Blue} reference, \, \sqbox{Green} $\varepsilon<0.1\%$:}
	\label{tab:eigNoFree}
\end{table}

\begin{table}[b]
	\centering	
	\begin{tabular}{|c|b g m|b g m|b g m|}
		\hline 
		\multirow{2}{*}{$n$}&   \multicolumn{3}{c|}{CSFS}	& \multicolumn{3}{c|}{SSFS} &	\multicolumn{3}{c|}{FSFS}	\\				
		&\cellcolor{white}Leissa & \cellcolor{white}Hessian & \cellcolor{white}divDiv  &\cellcolor{white}Leissa & \cellcolor{white}Hessian & \cellcolor{white}divDiv  &\cellcolor{white}Leissa & \cellcolor{white}Hessian & \cellcolor{white}divDiv \\
		\hline
		$\widehat{\omega}_{1}$&	12.69&	12.681&	\cellcolor{Orange}14.336&	11.68&	11.679&	\cellcolor{Orange}13.136&	9.631&	9.630&	\cellcolor{Orange}11.110\\
		$\widehat{\omega}_{2}$&	33.06&	33.041&	33.895&	27.76&	27.732&	28.681&	16.13&	16.117& \cellcolor{Orange}18.421\\
		$\widehat{\omega}_{3}$&	41.7 &	41.692&	43.791&	41.2 &	41.185&	43.073&	36.72&	36.683&	38.058\\
		$\widehat{\omega}_{4}$&	63.01&	62.982&	64.753&	59.07&	59.027&	60.764&	38.94&	38.939&	40.808\\
		$\widehat{\omega}_{5}$&	72.4 &	72.371&	\cellcolor{GreenYellow}72.756&	61.86&	61.825&	\cellcolor{GreenYellow}62.298&	46.74&	46.709&	\cellcolor{Orange}50.252\\
		$\widehat{\omega}_{6}$&	90.61&	90.602&	92.960&	90.29&	90.283&	92.436&	70.75&	70.666&	73.550\\	
		\hline 
	\end{tabular} 	
	\caption{Eigenvalues obtained with 5 Bell element per side for $\nu=0.3$, considering either the $\mathrm{Grad} \circ \mathrm{grad}$ formulation \eqref{eq:WF_Kir_Dyn}, either the $\mathrm{Div} \circ \mathrm{div}$ formulation \eqref{eq:WF_Kir_Kin}. For comparison reference \cite{LeissaRect} is considered: \\
		\sqbox{Blue} reference, \, \sqbox{Green} $\varepsilon<0.1\%$ , \, \sqbox{GreenYellow} $\varepsilon<1\%$, \, 			\sqbox{Yellow} $\varepsilon<5\%$, \, \sqbox{Orange} $\varepsilon<15\%$.}
	\label{tab:eigFree}
\end{table}


\subsection{Time-domain Simulations}
\begin{figure}[t]%
\centering
\subfloat[][$w(t = 0.25 \, t_{\text{fin}})$]{
	\label{fig:sim1-a}	\includegraphics[width=0.45\linewidth]{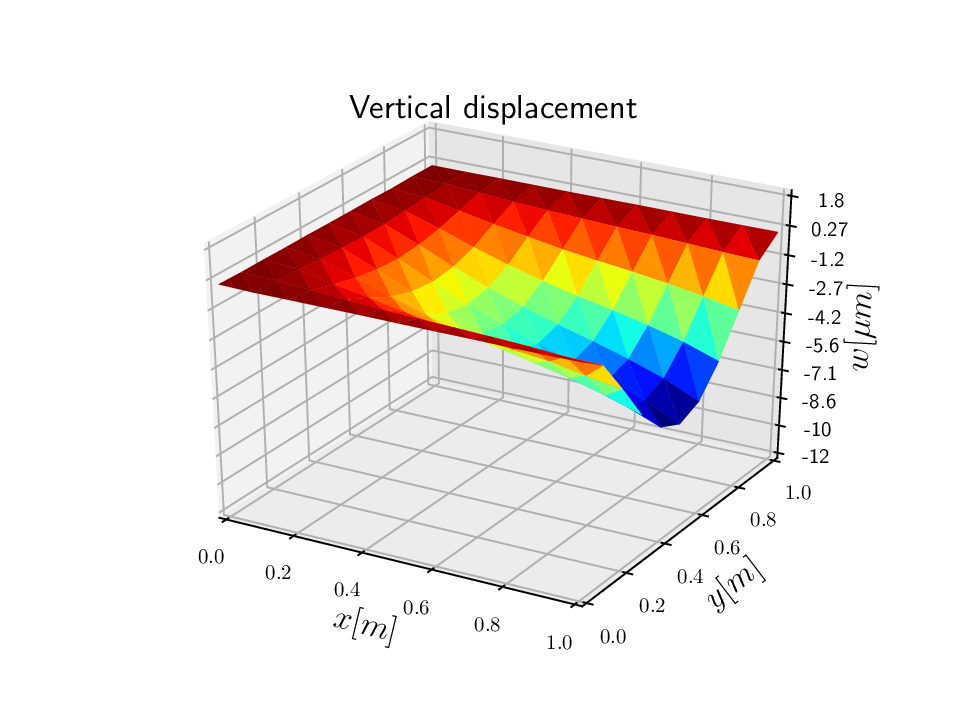}}
\hspace{8pt}
\subfloat[][$w(t = 0.50 \, t_{\text{fin}})$]{
	\label{fig:sim1-b}	\includegraphics[width=0.45\linewidth]{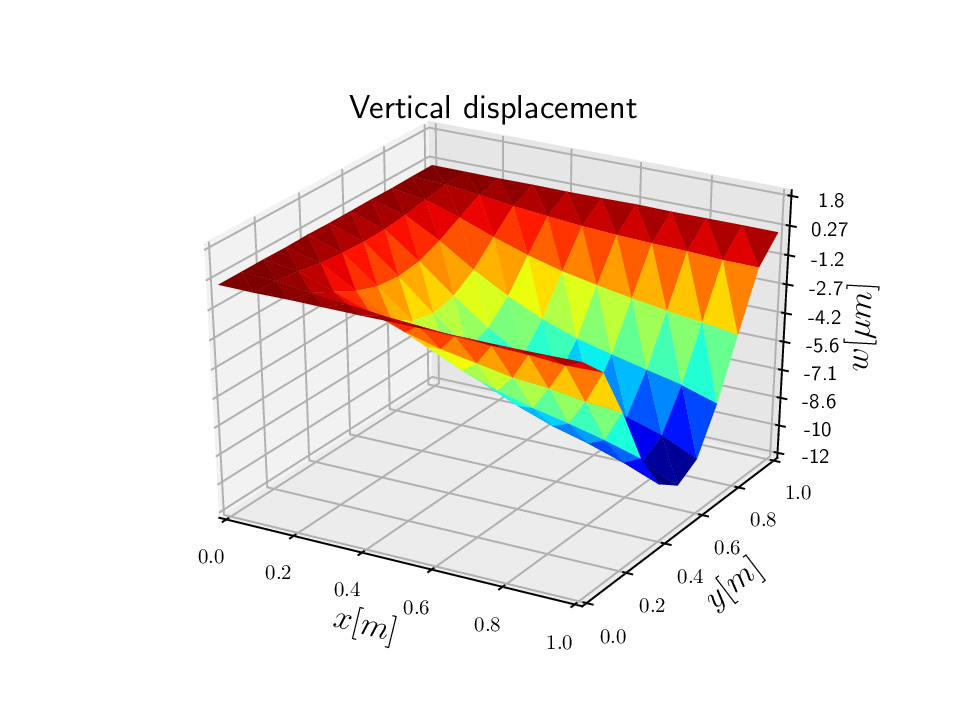}} \\
\subfloat[][$w(t = 0.75 \, t_{\text{fin}})$]{
	\label{fig:sim1-c}	\includegraphics[width=0.45\linewidth]{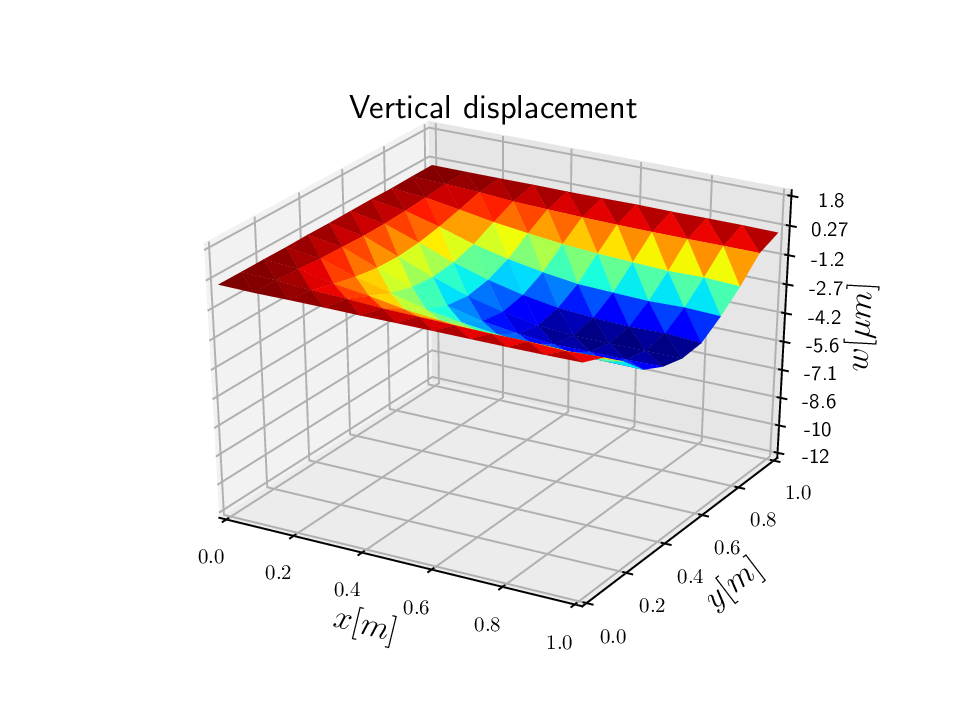}}
\hspace{8pt}
\subfloat[][$w(t = t_{\text{fin}})$]{
	\label{fig:sim1-d} \includegraphics[width=0.45\textwidth]{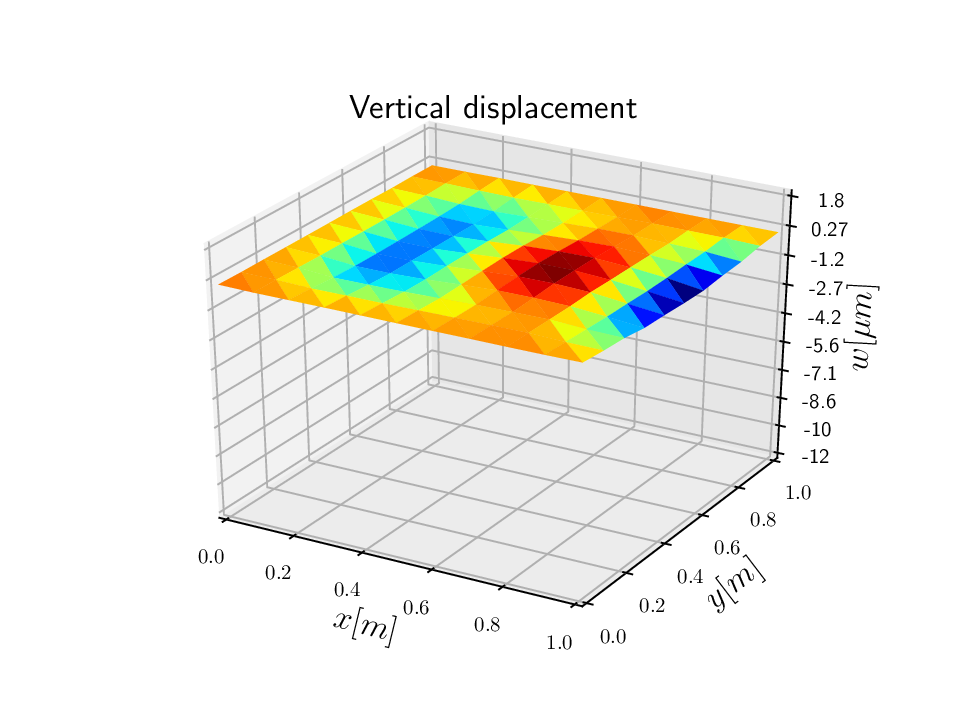}} 
\caption[Snapshots of the displacement field]{Snapshots for Simulation $n^\circ 1$.}%
\label{fig:sim1}%
\end{figure}

\begin{figure}[t]%
\centering
\subfloat[][$w(t = 0.25 \, t_{\text{fin}})$]{
	\label{fig:sim2-a} 	\includegraphics[width=0.45\textwidth]{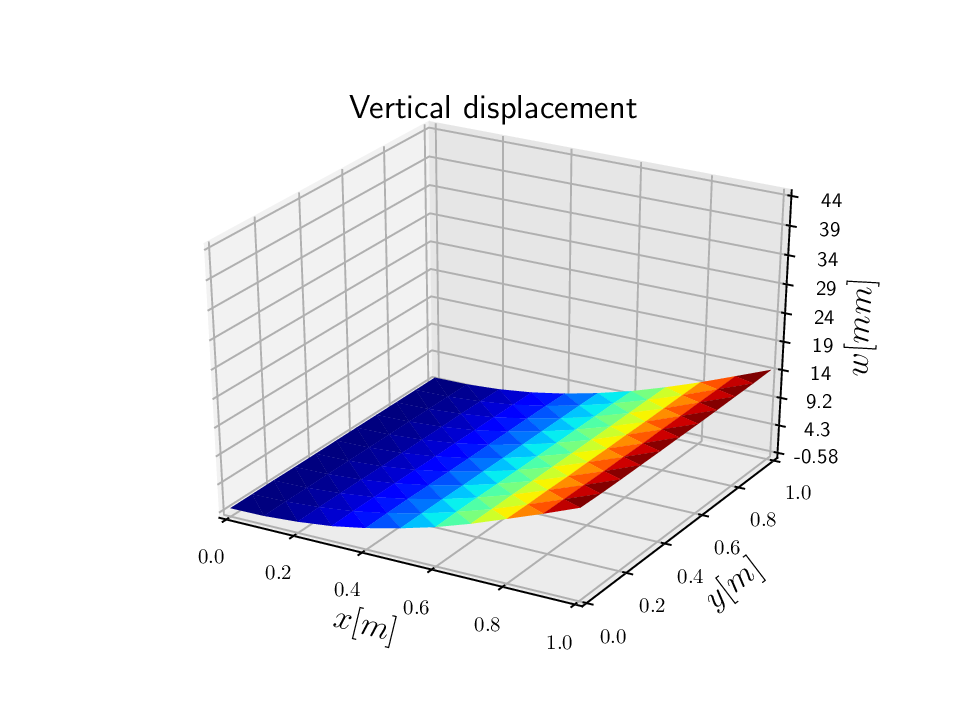}}
\hspace{8pt}
\subfloat[][$w(t = 0.50 \, t_{\text{fin}})$]{
	\label{fig:sim2-b} 	\includegraphics[width=0.45\textwidth]{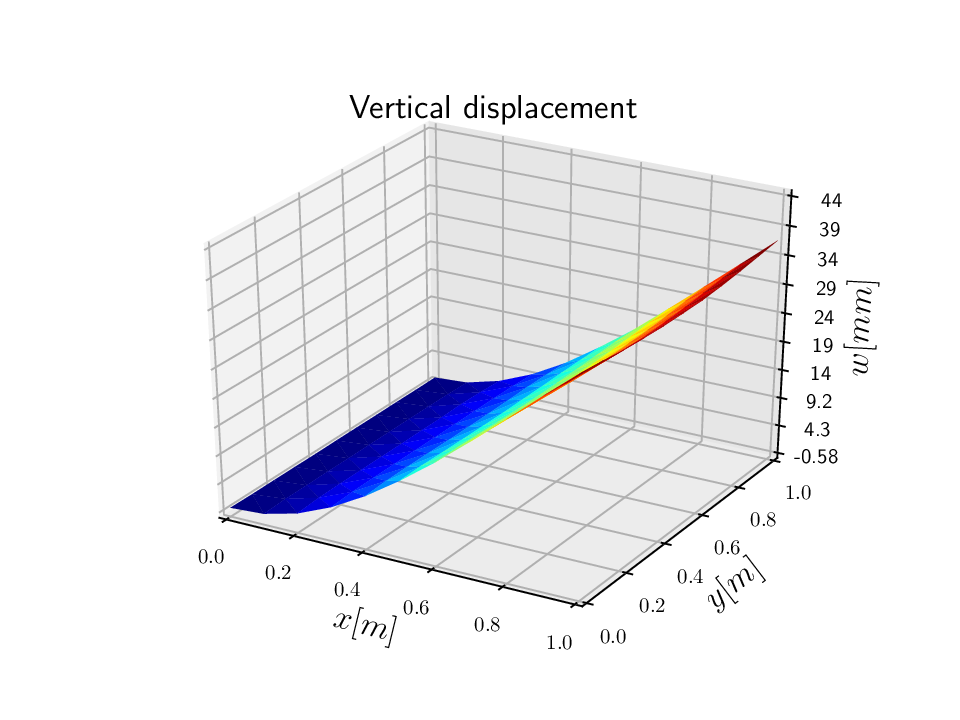}} \\
\subfloat[][$w(t = 0.75 \, t_{\text{fin}})$]{
	\label{fig:sim2-c}  \includegraphics[width=0.45\textwidth]{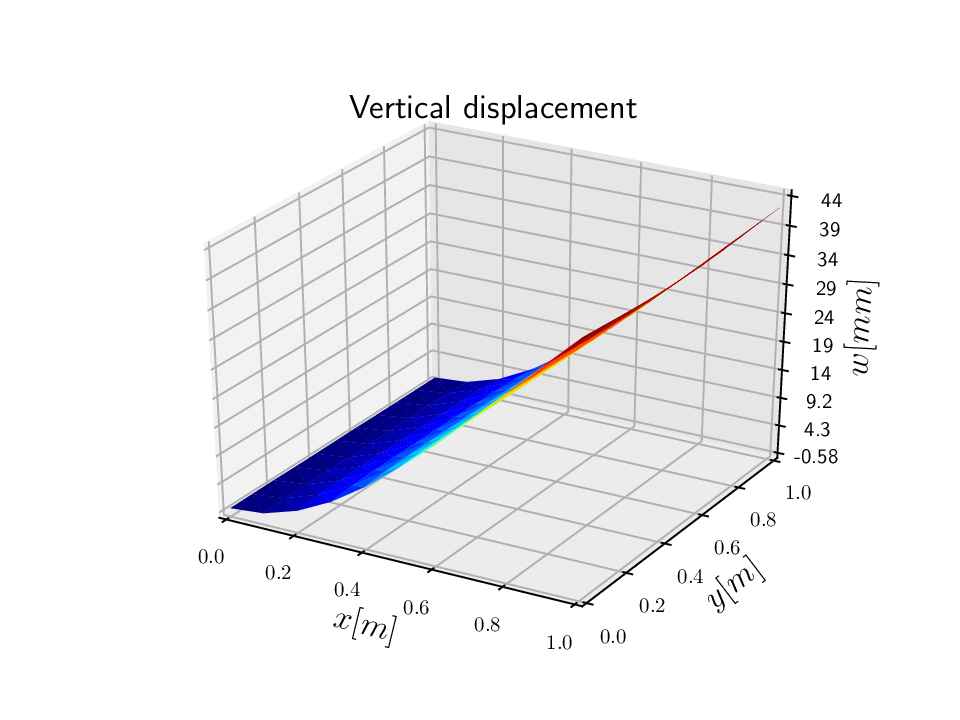}}
\hspace{8pt}
\subfloat[][$w(t = t_{\text{fin}})$]{
	\label{fig:sim2-d}	\includegraphics[width=0.45\textwidth]{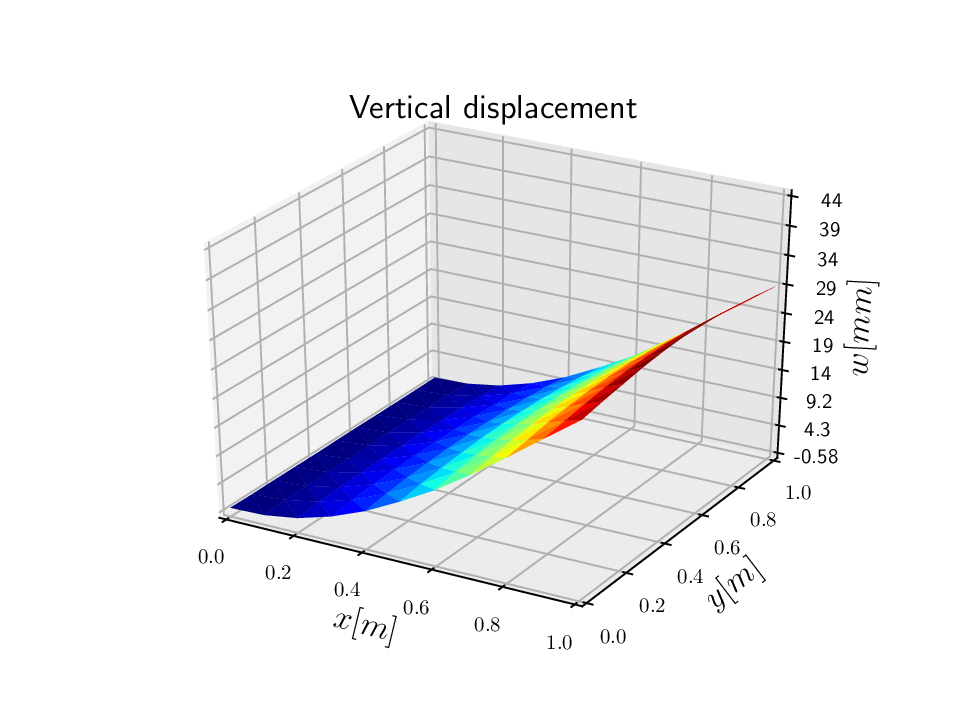}} 
\caption[Snapshots of the displacement field]{Snapshots for Simulation $n^\circ 2$.}%
\label{fig:sim2}%
\end{figure}

\begin{table}[t]
	\centering
	\begin{tabular}{|c|c|}
		\hline 
		\multicolumn{2}{|c|}{Plate Parameters} \\
		\hline 
		$E$ & $70\; [GPa]$ \\ 
		$\rho$ & $2700\; [kg/m^3]$ \\
		$\nu$& 0.35 \\ 
		$h/L$& 0.05 \\ 
		$L$& $1\; [m]$ \\ 
		\hline 
	\end{tabular} 
	\begin{tabular}{|c|c|}
		\hline 
		\multicolumn{2}{|c|}{Simulation Parameters} \\  
		\hline 
		Integrator & St\"ormer-Verlet \\ 
		$\Delta t$ & $0.001\;  [ms]$ \\ 
		$t_{\text{end}}$& $10 \; [ms]$ \\ 
		N$^\circ$ Elements & 5 \\
		FE space & $H_{r = L/5}^2(\mathbb{P}_5, \Omega) \text{ for }\bm{e} \times H_{r = L/5}^1(\mathbb{P}_2, \partial\Omega)\text{ for }\bm{\lambda}$ \\
		\hline 
	\end{tabular} 
	\captionsetup{width=0.95\linewidth}
	\vspace{1mm}
	\captionof{table}{Physical parameters and simulations settings.}
	\label{tab:par}
\end{table}
\begin{figure}[t]%
	\centering
	\subfloat[][Simulation $n^\circ 1$]{
		\label{fig:sim1-H}
		\includegraphics[width=0.47\textwidth]{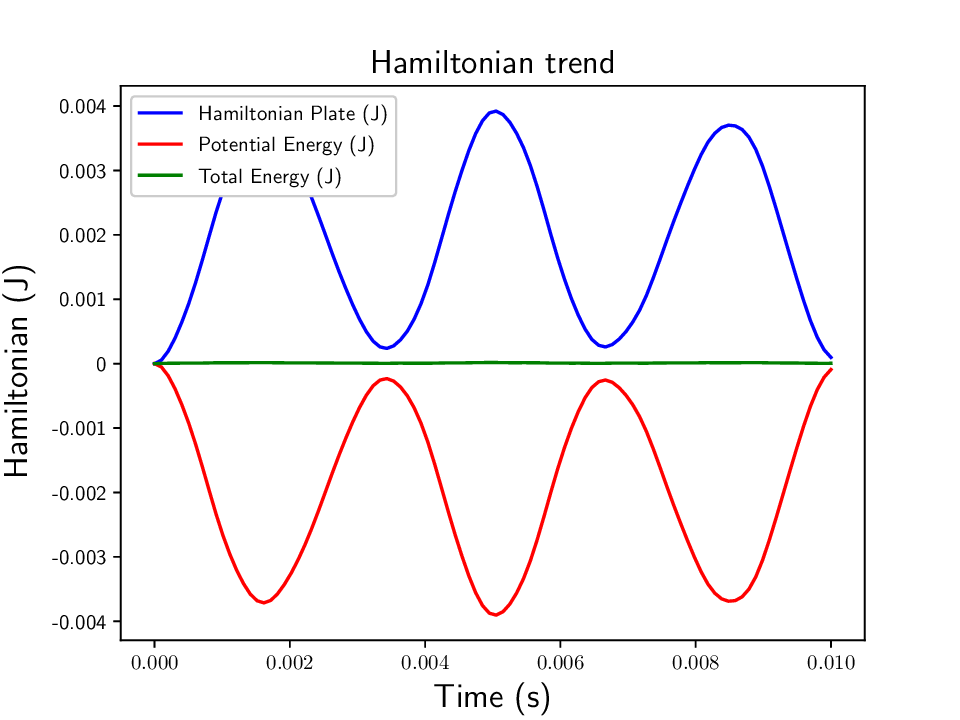}}
	\hspace{8pt}
	\subfloat[][Simulation $n^\circ 2$]{
		\label{fig:sim2-H}
		\includegraphics[width=0.45\textwidth]{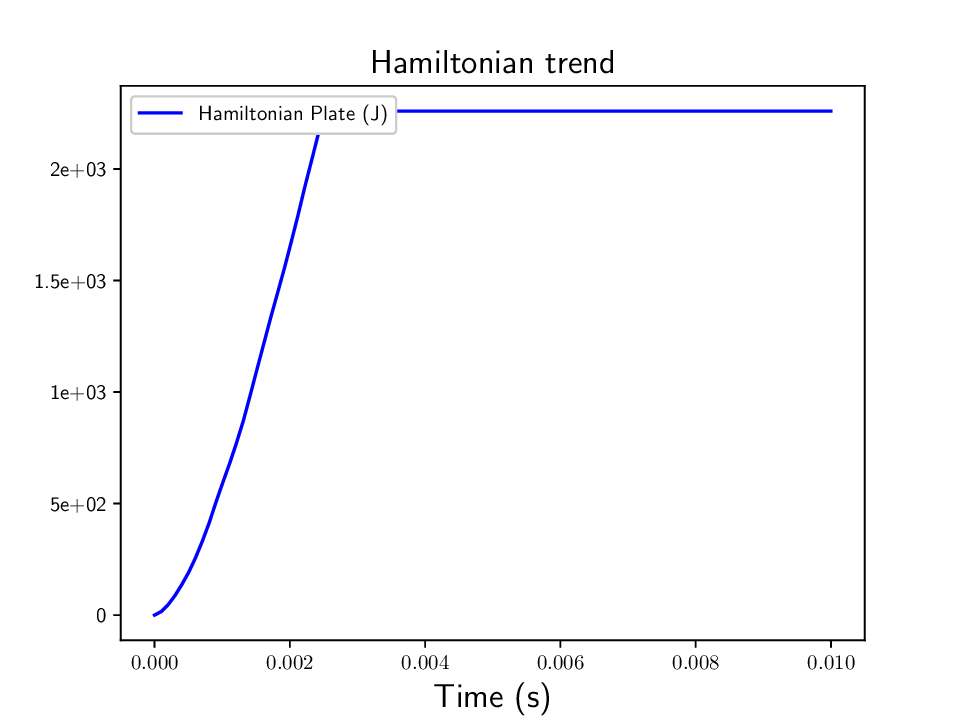}} \\
	\caption[Hamiltonian]{Hamiltonian trend for the two simulations.}
	\label{fig:Hamiltonian}
\end{figure}
In this analysis we consider a square plate, subject either to a non null shear force on the boundaries either to a distributed force over the domain. The physical parameters and simulation settings are reported in Table \ref{tab:par}. The energy variables and Lagrange multipliers are discretized using Bell shape functions (regular mesh of five elements for each side) and second order Lagrange polynomials respectively. The St\"ormer-Verlet time integrator is employed, so that the symplectic structure is preserved. Two different simulations with different boundary conditions are considered. The initial conditions are set to zero for each variable. For the first simulation a plate subject to gravity is considered. 
For this simulation, the following  boundary conditions, corresponding to the case CCCF are considered:
\begin{equation} \text{Simulation $n^\circ \, 1$} \;
\begin{cases}
w_t = 0, \, \diffp{w_t}{n} = 0,  \quad &\text{for } x=0, y=0 \text{ and } y=1 \\
\widetilde{q}_n = 0, \, M_{nn}= 0, \quad &\text{for } x=1\\
\end{cases}
\end{equation}
Since the solicitation admits a potential the Hamiltonian does not represent the total energy, that now includes the potential energy, whose expression is given by:
\begin{equation}
E_p = \int_{\Omega} \rho h g w \d{\Omega}, 
\end{equation}
where $w$ is the vertical displacement field and $g = 10 \, [m/s^2]$ is the gravity acceleration. \\
For the second simulation the following  boundary conditions are considered:
\begin{equation} \text{Simulation $n^\circ \, 2$} \;
\begin{cases}
w_t = 0, \, \diffp{w_t}{n} = 0, \quad &\text{for } x=0,\\
q_n = +f(1, t), \, M_{nn} = M_{ns}= 0, \quad &\text{for } x=1,\\
q_n = +f(x, t), \, M_{nn} = M_{ns}= 0, \quad &\text{for } y=0,\\
q_n = +f(x, t), \, M_{nn} = M_{ns}= 0, \quad &\text{for } y=1,\\
\end{cases}
\end{equation}
where the excitation $f(x, t)$ is computed as:
\begin{equation}
f(x, t) = \begin{cases}
10^5 \, x \; [Pa \cdot m], \; &\forall t< 0.25 \, t_{\text{end}}, \\
0, \; &\forall t\geq 0.25 \, t_{\text{end}}. \\
\end{cases}
\end{equation} 
In this case inhomogeneous boundary conditions are considered. Snapshots of the vertical displacement are reported in Figs. \ref{fig:sim1}, \ref{fig:sim2}. This field is obtained from the velocity field $e_w = \diffp{w}{t}$ by applying the trapezoidal rule integration. For both simulations, the output is consistent with the imposed BC and with the physical intuition of the observed phenomenon. The symplectic integration has been used to demonstrate numerically the conservation of total energy, as it can be noticed in Fig. \ref{fig:Hamiltonian}. 

\section*{Conclusions and Future Perspectives}
In this paper the port-Hamiltonian formulation of the Kirchhoff plate was detailed with the equivalent vectorial and a tensorial representation. The tensorial formalism allowed showing the adjointness relation between the double divergence of a symmetric tensor and the Hessian of a scalar field. This result represent a appealing novelty for mathematical working on functional a analysis. Moreover, many features of the PFEM are of interest:
\begin{itemize}
	\item its capability of preserving the port-Hamiltonian structure;
	\item the natural derivation of boundary port variables as inputs;
	\item the possibility of dealing with mixed boundary conditions inside the framework of  port-Hamiltonian descriptor systems PHDAEs  detailed in \cite{beattie2018linear};
	\item the easy implementability of the method using standard Finite Element libraries (Firedrake \cite{firedrake} in this case);
\end{itemize}   
The computation of eigenvalues with different boundary conditions and the numerical simulations demonstrate the validity of the proposed model. \\ \\
The model presented in this paper should be completed with a precise analysis of the well-posedness, in the input-output sense. It must be appointed that a complication arises in this formulation. The differential operator of the PH model of the Kirchhoff plate requires the momenta to belong to the space $H^{\text{div Div}}(\Omega, \mathbb{R}^{d \times d}_{\text{sym}})$. To the best of our knowledge this space was never analyzed in the mathematical literature and a precise study of its peculiarities is needed. The results obtained in \cite{waveEqZwart} for the wave equation in $\mathbb{R}^d$ could be generalize to the second order differential operator presented herein. A numerical analysis focusing on the convergence of appropriate finite elements should be carried out. \\ \\
The discretization procedure details in the paper open new scenarios on the interconnection of PH systems. Starting from the results stated in \cite{ShaftIntInfinite}, this system may be interconnected over its boundary to other finite or infinite dimensional PH systems, such as rigid bodies or other flexible appendages. This may find useful applications in simulating a multi-body environment for spatial applications, like the attitude motion of a satellite with flexible solar panels \cite{aoues:hal-01738092}. Since no causality is imposed on the boundary ports, the discretization method herein proposed allows the construction of arbitrarily complex connections among different modules. This feature is particularly appealing for complex applications.

\section*{Acknowledgments}
This work is supported by the project ANR-16-CE92-0028, entitled {\em Interconnected Infinite-Dimensional systems for Heterogeneous Media}, INFIDHEM, financed by the French National Research Agency (ANR) and the Deutsche Forschungsgemeinschaft (DFG). Further information is available at \\ {\url{https://websites.isae-supaero.fr/infidhem/the-project}}. \\
Moreover the author would like to thank Michel Sala\"un for the fruitful and insightful discussions

\bibliographystyle{unsrt}
\bibliography{biblio_Kirchh_Revision}

\end{document}